# The Feynman graph representation of convolution semigroups and its applications to Lévy statistics

HANNO GOTTSCHALK[1,*], BOUBAKER SMII[1,**] and HORST THALER[2]

[1] *Institut für angewandte Mathematik, Wegelerstr. 6, D-51373 Bonn, Germany. E-mail: *gottscha@wiener.iam.uni-bonn.de; **boubaker@wiener.iam.uni-bonn.de*
[2] *Dipartimento di Matematica e Informatica, Via Madonna delle Carceri 9, I 62032 Camerino, Italy. E-mail: horst.thaler@unicam.it*

We consider the Cauchy problem for a pseudo-differential operator which has a translation-invariant and analytic symbol. For a certain set of initial conditions, a formal solution is obtained by a perturbative expansion. The series so obtained can be re-expressed in terms of generalized Feynman graphs and Feynman rules. The logarithm of the solution can then be represented by a series containing only the connected Feynman graphs. Under some conditions, it is shown that the formal solution uniquely determines the real solution by means of Borel transforms. The formalism is then applied to probabilistic Lévy distributions. Here, the Gaussian part of such a distribution is re-interpreted as a initial condition and a large diffusion expansion for Lévy densities is obtained. It is outlined how this expansion can be used in statistical problems that involve Lévy distributions.

*Keywords:* Borel summability; convolution semigroups; Feynman graphs and rules; Lévy distributions; maximum likelihood principle

## 1. Introduction and overview

It is a well-known problem that there is no closed formula for the transition densities of general Lévy processes. In fact, though there are many examples of Lévy laws where the density is known (see, e.g., [1, 9, 19]), a general formula, even in the sense of series expansions, seems to be missing. Small-time expansions for pure jump Lévy densities have been considered (e.g., in [2, 17]). But such expansions require the knowledge of convolution powers of the Lévy measures which, in most cases, will be hard to calculate analytically. Expansions are available for stable laws (cf. [19, 21, 22] and references therein).

It seems that there is no expansion formula for Lévy densities in the case where the Lévy process has moments of all orders that is simple enough to be applicable in statistics.

In this work, we provide such a formula in the sense of Borel summable asymptotic expansions in two cases. Either the Lévy process has to be started in a initial distribution that is given by a density of special form or it can be started in a point but has to have a (large) diffusive part. The techniques employed stem from the renormalization group of statistical mechanics [18] and the recent proposal of a generalized Feynman graph calculus [5, 6]. Some of our ideas are similar to recent work of physicists motivated by Feynman path integration (cf. [9], Chapter 20.4.6, for a large diffusion expansion for one-dimensional and symmetric Lévy densities). In [11], a different line of approach for obtaining asymptotics of densities given as solutions of pseudo-differential operators is presented. Furthermore, the expansion that we present is related, but not identical, to the well-known bootstrap expansion in statistics (cf., e.g., [3, 8] for an introduction to this flourishing field). However, it seems that non-symmetric densities on $\mathbb{R}^d$, Borel summability, the graphic representation and so-called linked cluster theorems have not been considered before.

Most of our considerations are not necessarily restricted to Lévy processes, but hold for fairly general convolution semigroups. We study the Cauchy problem of the partial (pseudo-) differential equation

$$\begin{aligned} \frac{\partial}{\partial t}\Phi_t(\phi) &= \Psi(\nabla)\Phi_t(\phi), \qquad (t,\phi)\in\,]0,\infty[\,\times\mathbb{R}^d, \\ \Phi_0(\phi) &= f(\phi), \qquad \phi\in\mathbb{R}^d \end{aligned} \tag{1}$$

for a given initial state $f:\mathbb{R}^d\to\mathbb{C}$ and a (pseudo-) differential operator with constant coefficients $\Psi(\nabla)$.

We consider the case where $\Psi$ is the symbol of a (pseudo-) differential operator, that is, has an expansion

$$\Psi(\xi) = \sum_{n=0}^{\infty}\langle C^{(n)},\xi^{\otimes n}\rangle_n = \sum_{n=0}^{\infty}\sum_{X_1,\dots,X_n=1}^{d} C^{(n)}_{X_1\cdots X_n}\xi_{X_1}\cdots\xi_{X_n}, \tag{2}$$

$\xi=(\xi_1,\dots,\xi_d)\in\mathbb{R}^d$, that converges on some neighborhood of 0 and $f$ is of the form

$$f=\mathrm{e}^{-V}, \qquad V(\phi)=\sum_{p=0}^{\bar p}\langle\lambda^{(p)},\phi^{\otimes p}\rangle_p = \sum_{p=0}^{\bar p}\sum_{X_1,\dots,X_p=1}^{d}\lambda^{(p)}_{X_1\cdots X_n}\phi_{X_1}\cdots\phi_{X_p} \tag{3}$$

for $\phi=(\phi_1,\dots,\phi_d)\in\mathbb{R}^d$, where $\bar p$ is even and $\langle\lambda^{(\bar p)},\phi^{\otimes\bar p}\rangle_{\bar p}=\sum_{X_1,\dots,X_{\bar p}=1}^{d}\lambda^{(\bar p)}_{X_1\cdots X_{\bar p}}\phi_{X_1}\cdots\phi_{X_{\bar p}}>0$ for $\phi\neq 0$. We obtain a formal expansion of the solution of (1) in powers of $\lambda^{(p)}$ that can be represented in terms of generalized Feynman graphs [5, 6]. Here, $C^{(p)}$ ($C^{(p)}_{X_1\cdots X_p}$) and $\lambda^{(p)}$ ($\lambda^{(p)}_{X_1\cdots X_p}$) are symmetric tensors of $p$th degree over $\mathbb{R}^d$ (their respective components). $\langle\cdot,\cdot\rangle_p$ is the canonical scalar product on $(\mathbb{R}^d)^{\otimes p}$.

For example, in the case where $\lambda^{(p)} = 0$ unless $p = 2$ (this will be the most interesting case) and the symbol $\Psi(\xi)$ fulfills[1] $C^{(1)} = 0$, the lowest orders in this expansion are

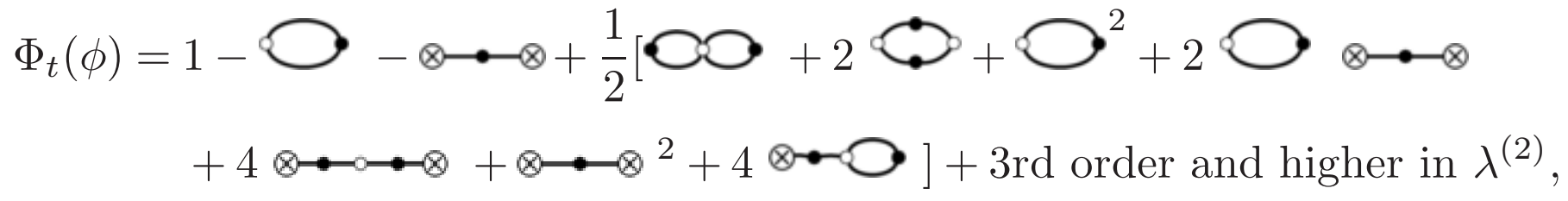

(4)

where each graph stands for a numerical expression that can be calculated in terms of the basic parameters $C^{(n)}_{X_1\cdots X_n}$, $\lambda^{(p)}_{X_1\cdots X_p}$, $\phi$ and $t$ by an easy algorithm ("Feynman rule") (cf. the following section).

This solution is the well-known graphical expansion for the renormalization group equation [18] for the case when $\Psi(\nabla)$ is a second-order partial differential operator (e.g., the Laplacian) and (1) is hence the heat equation. The recently discovered generalized Feynman graph calculus allows the generalization of this technique to a large class of symbols. This generalization is carried through in the present article, having in mind in particular Lévy processes with generators given by symbols of the form

$$\Psi(-\mathrm{i}\xi) = \mathrm{i}\langle a, \xi\rangle - \langle \xi, D\xi\rangle + z\int_{\mathbb{R}^d\setminus\{0\}} (\mathrm{e}^{\mathrm{i}\langle\varphi,\xi\rangle} - 1)\,\mathrm{d}r(\varphi), \tag{5}$$

which lead to equations (1) in the following jump-diffusion form (Lévy–Itô decomposition):

$$\frac{\partial\Phi_t}{\partial t}(\phi) = -\sum_{X=1}^{d} a_X \frac{\partial}{\partial\phi_X}\Phi_t(\phi) + \sum_{X_1,X_2=1}^{d} D_{X_1X_2}\frac{\partial^2}{\partial\phi_{X_1}\partial\phi_{X_2}}\Phi_t(\phi)$$
$$+ z\int_{\mathbb{R}^d\setminus\{0\}} [\Phi_t(\phi+\varphi) - \Phi_t(\phi)]\,\mathrm{d}r(\varphi). \tag{6}$$

Here, $a \in \mathbb{R}^d$, $D$ is a real and positive semidefinite $d\times d$ matrix, $z \geq 0$ and $r$ is a probability measure on $\mathbb{R}^d \setminus \{0\}$ such that its Fourier transform is entire analytic. The well-known interpretation of $a$ is as the drift vector. $D$ determines the diffusion part, whereas $z$ and $r$ give the frequency and distribution of jumps, respectively.

The restriction of the set of initial conditions in (3) is unwanted in many cases. It can be circumvented by splitting the symbol (5) into a pure jump and a pure diffusion part. After some time, the transition kernel of the heat semigroup generated by the diffusive part alone is of course a Gaussian density, which, up to normalization, is of the form given in (3). This simple observation makes it possible to develop a perturbative calculus for densities of Lévy processes starting in a point in the limit of large diffusion. Apart from this assumption, the perturbative formalism developed in this article is fairly universal and we therefore hope that it can contribute to a better handling of Lévy distributions in statistics.

[1] The assumption concerning $\Psi(\xi)$ is without loss of generality since a non-zero $C^{(1)}$ can be compensated by a shift of $\phi$.

The article is organized as follows. In Section 2, we develop a perturbative formalism for the solution of equation (1). In Section 2.1, we prove the differentiability of the solution $\Phi_t$ in the parameters $\lambda^{(p)}$ that define the initial condition and conclude that the series so obtained is asymptotic. Section 2.2 gives an expansion of the perturbation series in the coefficients of $\Psi$ which in the case of Lévy processes are just the cumulants. The rearrangement in terms of generalized Feynman graphs is presented in 2.3. In the following Section 2.4, we obtain an expression for the logarithm of the solution $\Phi_t$ in terms of connected Feynman graphs. This does not only reduce the number of Feynman graphs that one has to calculate, but is also of interest in connection with applications to the maximum likelihood method. A brief consideration of Borel summability of the solution for the case that the logarithm of the initial condition, $V$, is quadratic concludes this section (Section 2.5).

In Section 3, we apply the general results to the special case of Lévy distributions. Section 3.1 presents the large diffusion expansion of Lévy densities. In Section 3.2, we illustrate our method by a second-order calculation of the (non-normalized) Lévy distributions using Padé resummation. The result is used for a calculation of densities and quantiles that we compare with the Monte Carlo simulation of a Lévy process with compound or pure Poisson jump part.

Some technical material that is presumably standard to some readers, but possibly not to others, is deferred to Appendices A.1–A.3.

# 2. Generalized Feynman graphs and the convolution semigroup

## 2.1. A solution in the sense of formal series

After the previous rough description, let us now fix the mathematical details. Let the pseudo-differential operator $\Psi$ fulfill the following requirements:

(P1) the associated symbol $\Psi(-\mathrm{i}\xi)$ is analytic for $\xi$ in a neighborhood of 0;
(P2) $\Re\Psi(-\mathrm{i}\xi) \leq c\log(1+\|\xi\|)$, where $\xi \in \mathbb{R}^d$ and $\Re$ denotes the real part;
(P3) $\forall\alpha \in \mathbb{N}_0^d$, $\exists N_\alpha \in \mathbb{N}, b_\alpha > 0$ such that $|D_\alpha\Psi(-\mathrm{i}\xi)| \leq b_\alpha(1+\|\xi\|)^{N_\alpha}$ $\forall\xi \in \mathbb{R}^d$.

Here, $D_\alpha h(\phi) = \frac{\partial^{|\alpha|}}{\partial\phi_1^{\alpha_1}\cdots\partial\phi_d^{\alpha_d}}h(\phi)$. We furthermore assume, without loss of generality, that $\Psi(0) = 0$. If this is not the case and $\Psi(0) = b$ (i.e., if the associated Lévy process undergoes killing or creation with rate $b$), then the solution $\Phi_t$ of equation (1) can be obtained from the solution $\tilde{\Phi}_t$ of the equation (1) with $\Psi$ replaced by $\tilde{\Psi} = \Psi - \Psi(0)$ and setting $\Phi_t(\phi) = \mathrm{e}^{t\Psi(0)}\tilde{\Phi}_t(\phi)$.

We take specific initial conditions that fulfill

(I1) $f(\varphi) = \mathrm{e}^{-V(\varphi)}$;
(I2) $V$ is a polynomial as given in equation (3).

Certainly, the restriction in the initial conditions by (I1)–(I2) is disturbing. Nevertheless, (I1)–(I2) are general enough to completely determine the convolution kernel $\nu_t$, as one can construct approximating sequences of the Dirac delta distribution in that class (e.g., take $V_\epsilon(\varphi) = \frac{\epsilon}{2}\|\varphi\|^2 - \frac{d}{2}\log(2\pi\epsilon)$, $\epsilon > 0$).

Let $\mathcal{F}: \mathcal{S}(\mathbb{R}^d) \to \mathcal{S}(\mathbb{R}^d)$ be the Fourier transform on the space of Schwartz test functions over $\mathbb{R}^d$, that is,

$$\mathcal{F}(f)(\xi) = \int_{\mathbb{R}^d} f(\varphi) \mathrm{e}^{\mathrm{i}\langle \xi, \varphi\rangle}\, \mathrm{d}\varphi, \qquad f \in \mathcal{S}(\mathbb{R}^d). \tag{7}$$

We denote the space of tempered distributions, that is, the topological dual of the Schwartz space $\mathcal{S}(\mathbb{R}^d)$, by $\mathcal{S}'(\mathbb{R}^d)$. By the duality given by the $L^2$ scalar product on $\mathbb{R}^d$ with respect to the Lebesgue measure, $\mathcal{S}(\mathbb{R}^d)$ is densely embedded in $\mathcal{S}'(\mathbb{R}^d)$ and $\mathcal{F}$ extends uniquely to the latter space by continuity with respect to the weak topology [4]. We denote this extension by the same symbol $\mathcal{F}$.

By condition (P2), $\mathrm{e}^{t\Psi(-\mathrm{i}\xi)}$ is polynomially bounded for all $t$, that is, $|\mathrm{e}^{t\Psi(-\mathrm{i}\xi)}| \leq (1 + \|\xi\|)^{ct}$, hence $\mathrm{e}^{t\Psi(-\mathrm{i}\xi)}$ is in $\mathcal{S}'(\mathbb{R}^d)$ for all $t \geq 0$ and $c > 0$ as in that condition. Furthermore, by conditions (I1) and (I2) – the positivity condition on $\lambda^{(\bar{p})}$ in particular – $f_\beta = \mathrm{e}^{-\beta V}$ is in the test function space $\mathcal{S}(\mathbb{R}^d)$ for $\beta > 0$. Hence, for $f = f_1$, the convolution $\Phi_t(\phi) = \nu_t * f(\phi) = \langle \nu_t, f_\phi\rangle$ is well defined, where $\nu_t \in \mathcal{S}'(\mathbb{R}^d)$ is the inverse Fourier transform of $\mathrm{e}^{t\Psi(-\mathrm{i}\xi)}$ and $f_\phi(\varphi) = f(\varphi - \phi)$.

**Lemma 2.1.** *$\Phi_t$ is the (unique) solution of the Cauchy problem (1).*

**Proof.** The proof is a standard argument based on the Fourier transform of tempered distributions. For the convenience of the reader, the details are given in Appendix A.1. □

Following Poincaré, a formal power series $\sum_{m=0}^{\infty} \beta^m a_m$ is called an asymptotic series for the function $h(\beta)$, $\beta > 0$, at $\beta = 0$ if, for $N \in \mathbb{N}$,

$$\lim_{\beta \searrow 0} \beta^{-N} \left| h(\beta) - \sum_{m=0}^{N} \beta^m a_m \right| = 0. \tag{8}$$

Here, we want to expand the solution $\Phi_t$ in powers of $V$. In many cases, this expansion is not convergent. However, we show that such an expansion is asymptotic in the above sense. To this aim, let, for $\beta > 0$, $\Phi_t^\beta = \nu_t * f_\beta$. Expansion in powers of $V$ then means to expand in the auxiliary parameter $\beta$ and then set $\beta = 1$. By definition, the series in $V$ is asymptotic to $\Phi_t$ if and only if the $\beta$ expansion exists and is asymptotic to $\Phi_t^\beta$.

The *moments* of the distribution $\nu_t$ are defined as follows:

$$\langle \varphi_{X_1} \cdots \varphi_{X_n} \rangle_{\nu_t} = \frac{\partial^n}{\partial \xi_{X_1} \cdots \partial \xi_{X_n}} \mathrm{e}^{t\Psi(\xi)} \bigg|_{\xi = 0}, \qquad n \in \mathbb{N},\ X_1, \ldots, X_n \in \{1, \ldots, d\}. \tag{9}$$

For a polynomial $P$ in $\varphi$, the expression $\langle P(\varphi)\rangle_{\nu_t}$ is defined by linearity of the moments in the polynomial's coefficients.

**Proposition 2.2.** *The $V$-expansion of $\Phi_t$ is given by*

$$\Phi_t(\phi) = \sum_{m=0}^{\infty} \frac{(-1)^m}{m!} \langle V_\phi^m(\varphi)\rangle_{\nu_t}, \qquad V_\phi(\varphi) = V(\varphi - \phi). \tag{10}$$

*Equation (10) is to be understood in the sense that the right-hand side gives the asymptotic series for the left-hand side.*

**Proof.** The statement is evident in the sense of formal power series. That the right-hand side is asymptotic follows from the fact that $\Phi_t^\beta(\phi)$ is in $C^\infty((0,\infty))$ with respect to $\beta$ and for any $m \in \mathbb{N}$, $\frac{\mathrm{d}^m}{\mathrm{d}\beta^m}\Phi_t^\beta(\phi)$ can be continuously extended to $\beta = 0$ such that $\lim_{\beta\searrow 0} \frac{\mathrm{d}^m}{\mathrm{d}\beta^m}\Phi_t^\beta(\phi) = (-1)^m\langle V_\phi^m(\varphi)\rangle_{\nu_t}$ (cf. Appendix A.1 for the details). A posteriori, we conclude that the extension of $\Phi_t^\beta$ is $C^\infty([0,\infty))$ in $\beta$ and apply Taylor's lemma to conclude the argument. □

In general, the asymptotic expansion is not convergent. This applies in particular to $\Psi(-\mathrm{i}\xi)$ given in (5), where $\Psi$ is the generator of a Lévy process. We will return to this point in Section 2.5, where Borel summability is proved in an important case. By the fact that the expansion is asymptotic, it is already clear that the first $N$ terms in (10) give a good approximation to the actual solution for $V$ small and $\phi$ not too large, an error estimate being provided by Taylor's formula. In order to extract reliable data for larger $V$, resummation techniques have to be applied. Some first steps in this direction are given in Section 3.2.

Obviously, all partial differential operators of the form $\Psi(\nabla) = \sum_{n=0}^{\bar{n}} \sum_{X_1,\ldots,X_n=1}^{d} C^{(n)}_{X_1\cdots X_n} \times \frac{\partial^n}{\partial\varphi_{X_1}\cdots\partial\varphi_{X_n}}$ of even order $\bar{n}$ with $(-1)^{\bar{n}/2}\langle C^{(\bar{n})}, \psi^{\otimes\bar{n}}\rangle_{\bar{n}} > 0$ for $\xi \neq 0$ or of odd order $\bar{n}$ with $C^{(\bar{n})}$ real and $(-1)^{(\bar{n}-1)/2}\langle C^{(\bar{n}-1)}, \psi^{\otimes\bar{n}-1}\rangle_{\bar{n}-1} > 0$ for $\xi \neq 0$ fulfill the conditions (P1)–(P3). This, of course, includes the heat equation ($\bar{n} = 2$). More interestingly, not only local generators can be treated, but also generators of Lévy processes that have a symbol that is analytic at 0.

**Proposition 2.3.** *The generator (5) with $r$ a probability measure with analytic Fourier transform at 0 fulfills the conditions* (P1)–(P3).

**Proof.** (P1), analyticity at 0, is immediate from the conditions. The real part of (5) is obviously bounded from above by $2z$, which establishes (P2). To see (P3), let us note that the partial derivatives of the third part in (5) are of the form $z\int_{\mathbb{R}^d\setminus\{0\}} \varphi^\alpha \mathrm{e}^{\mathrm{i}\langle\xi,\varphi\rangle}\,\mathrm{d}r(\varphi)$. As $r$ has analytic Fourier transform at zero, all moments of $r$ exist. This implies differentiability everywhere and even that $D_\alpha\Psi(-\mathrm{i}\xi)$ is uniformly bounded in $\xi \in \mathbb{R}^d$. □

For the sake of mathematical completeness, we also give a sufficient condition for the convergence of the expansion (10). Its practical value is limited, however, as it excludes the most interesting examples.

**Proposition 2.4.** *If the symbol* $\Psi$ *does not necessarily fulfill* (P3) *but fulfills* (P3′)*:* $\Psi(-\mathrm{i}\xi)$ *is entire analytic and* $\Psi(-\mathrm{i}\xi)\le K(1+\|\Im\xi\|_{\max})$*,* $\xi\in\mathbb{C}^d$*, for some* $K>0$*, then the expansion* (10) *is convergent.*

**Proof.** The condition implies that $\nu_t$ has an entire analytic Fourier transform $|\mathrm{e}^{t\Psi(\xi)}|\le \mathrm{e}^{Kt}\mathrm{e}^{tK\|\Im\psi\|_{\max}}$. Application of the Paley–Wiener–Schwartz theorem [14] implies that $\nu_t$ has support in a square of length $2tK$ centered at the origin and hence compact support. Let $\chi$ be a test function with compact support and $\chi|_{\mathrm{supp}\,\nu_t}=1$. Then $\Phi_t^\beta(\phi)=\nu_t * f_\beta(\phi)=\langle\nu_t,\chi\cdot\mathrm{e}^{-\beta V_\phi}\rangle$. As the series expansions in $\beta$, $\beta\in\mathbb{C}$, of $\mathrm{e}^{-\beta V_\phi}$ and all of its derivatives converge compactly on $\mathrm{supp}\,\chi$, the mapping $\mathbb{C}\ni\beta\to\chi\cdot\mathrm{e}^{-\beta V_\phi}\in\mathcal{S}(\mathbb{R}^d)$ is entire analytic in $\beta$. Since $\nu_t\in\mathcal{S}'(\mathbb{R}^d)$, the same applies to $\Phi_t^\beta(\phi)$. □

## 2.2. Expansion in the operator's coefficients

In this section, we define the truncated moments of the convolution kernel $\nu_t$ as partial derivatives of $t\Psi$ – as the coefficients of $\Psi$ times $t$, in other words. A result known as the *linked cluster theorem* (sometimes also named the *cumulant expansion*) establishes a combinatorial connection between the moments $\langle\varphi_{X_1}\cdots\varphi_{X_n}\rangle_{\nu_t}$ and the truncated moments $\langle\varphi_{X_1}\cdots\varphi_{X_n}\rangle_{\nu_t}^T$. The asymptotic expansion of the solution of the Cauchy problem (1) can thus be directly expressed in terms of the coefficients of the symbol $C^{(n)}_{X_1\cdots X_n}$, the time $t$ and the coefficients $\lambda^{(p)}_{X_1\cdots X_p}$ of the logarithm of the initial condition.

***Definition 2.5.*** *The truncated moments of the distribution* $\nu_t$ *are defined as*

$$\langle\varphi_{X_1}\cdots\varphi_{X_n}\rangle_{\nu_t}^T=t\frac{\partial^n}{\partial\xi_{X_1}\cdots\partial\xi_{X_n}}\Psi(\xi)|_{\xi=0},\qquad n\in\mathbb{N},\ X_1,\ldots,X_n\in\{1,\ldots,d\}. \tag{11}$$

By (2), one obtains that the truncated moments are the coefficients of the (pseudo-) differential operator $\Psi$, that is,

$$\langle\varphi_{X_1}\cdots\varphi_{X_n}\rangle_{\nu_t}^T=tC^{(n)}_{X_1\cdots X_n}. \tag{12}$$

The following classical theorem gives the relation between the ordinary moments $\langle\varphi_{X_1}\cdots\varphi_{X_n}\rangle_{\nu_t}$ and the truncated ones introduced above.

**Theorem 2.6 (Linked cluster).** *For* $n\in\mathbb{N}$*, we have*

$$\langle\varphi_{X_1}\cdots\varphi_{X_n}\rangle_{\nu_t}=\sum_{\substack{I\in\mathcal{P}(\{1,\ldots,n\})\\ I=\{I_1,\ldots,I_k\}}}\prod_{l=1}^k\left\langle\prod_{j\in I_k}\varphi_{X_j}\right\rangle_{\nu_t}^T, \tag{13}$$

*where, for a finite set* $A$*,* $\mathcal{P}(A)$ *stands for the set of all partitions* $I$ *of* $A$ *into non-empty disjoint subsets* $\{I_1,\ldots,I_k\}$*.*

For the proof, we refer to the literature (e.g., [5, 16]).

**Proposition 2.7.** *The solution of the Cauchy problem $\Phi_t(\phi)$, $\phi \in \mathbb{R}^d$, has the following asymptotic expansion:*

$$\Phi_t(\phi) = \sum_{m=0}^{\infty} \frac{(-1)^m}{m!} \sum_{p_1,\dots,p_m=0}^{\bar{p}} \sum_{\substack{K \subseteq \Omega(p_1,\dots,p_m) \\ I \in \mathcal{P}(\Omega(p_1,\dots,p_m)\setminus K) \\ I=\{I_1,\dots,I_k\}}} \sum_{X_1^1,\dots,X_{p_1}^1=1}^{d} \cdots \sum_{X_1^m,\dots,X_{p_m}^m=1}^{d} \lambda^{(p_1)}_{X_1^1,\dots,X_{p_1}^1} \cdots \lambda^{(p_m)}_{X_1^m,\dots,X_{p_m}^m} \times \prod_{(s,q)\in K} (-\phi_{X_s^q}) \prod_{l=1}^{k} tC_{I_l}. \tag{14}$$

*Here,* $\Omega(p_1,\dots,p_m) = \bigcup_{n=1}^{m} \{(1,n),\dots,(p_n,n)\}$ *and* $C_{I_l} = C^{(q)}_{X_{a_1}^{b_1}\cdots X_{a_q}^{b_q}}$ *for* $I_l = \{(a_1,b_1),\dots,(a_q,b_q)\}$.

**Proof.** The proof is by a straightforward calculation using (11) and (12):

$$\begin{aligned}
\Phi_t(\phi) &= \sum_{m=0}^{\infty} \frac{(-1)^m}{m!} \langle V_\phi^m(\varphi)\rangle_{\nu_t} \\
&= \sum_{m=0}^{\infty} \frac{(-1)^m}{m!} \sum_{p_1,\dots,p_m=0}^{\bar{p}} \sum_{X_1^1,\dots,X_{p_1}^1=1}^{d} \cdots \sum_{X_1^m,\dots,X_{p_m}^m=1}^{d} \lambda^{(p_1)}_{X_1^1\cdots X_{p_1}^1} \cdots \lambda^{(p_m)}_{X_1^m\cdots X_{p_m}^m} \\
&\quad \times \langle (\varphi_{X_1^1} - \phi_{X_1^1}) \cdots (\varphi_{X_{p_1}^1} - \phi_{X_{p_1}^1}) \cdots (\varphi_{X_1^m} - \phi_{X_1^m}) \cdots (\varphi_{X_{p_m}^m} - \phi_{X_{p_m}^m}) \rangle_{\nu_t} \\
&= \sum_{m=0}^{\infty} \frac{(-1)^m}{m!} \sum_{p_1,\dots,p_m=0}^{\bar{p}} \sum_{X_1^1,\dots,X_{p_1}^1=1}^{d} \cdots \sum_{X_1^m,\dots,X_{p_m}^m=1}^{d} \lambda^{(p_1)}_{X_1^1,\dots,X_{p_1}^1} \cdots \lambda^{(p_m)}_{X_1^m,\dots,X_{p_m}^m} \qquad (15) \\
&\quad \times \sum_{K \subset \Omega(p_1,\dots,p_m)} \prod_{(s,q)\in K} (-\phi_{X_s^q}) \left\langle \prod_{(r,t)\in\Omega(p_1,\dots,p_m)\setminus K} \varphi_{X_t^r} \right\rangle_{\nu_t} \\
&= \sum_{m=0}^{\infty} \frac{(-1)^m}{m!} \sum_{p_1,\dots,p_m=0}^{\bar{p}} \sum_{X_1^1,\dots,X_{p_1}^1=1}^{d} \cdots \sum_{X_1^m,\dots,X_{p_m}^m=1}^{d} \lambda^{(p_1)}_{X_1^1\cdots X_{p_1}^1} \cdots \lambda^{(p_m)}_{X_1^m\cdots X_{p_m}^m} \\
&\quad \times \sum_{K \subset \Omega(p_1,\dots,p_m)} \prod_{(s,q)\in K} (-\phi_{X_s^q}) \sum_{\substack{I \in \mathcal{P}(\Omega(p_1,\dots,p_m)\setminus K) \\ I=\{I_1,\dots,I_k\}}} \prod_{l=1}^{k} tC_{I_l}.
\end{aligned}$$

□

Though the expansion is explicit, it is also rather messy. To simplify it, we introduce a graphical calculus in the next subsection.

If $d = 1$, the expression (14) simplifies considerably; one can drop all subscripts $X_1, \dots$ and thereby obtain the asymptotic expansion

$$\Phi_t(\phi) = \sum_{m=0}^{\infty} \frac{(-1)^m}{m!} \sum_{p_1,\dots,p_m=0}^{\bar{p}} \prod_{l=1}^{m} \lambda^{(p_l)} \sum_{k=0}^{p_1+\cdots+p_m} \binom{p_1+\cdots+p_m}{k} (-\phi)^k \\ \times \sum_{q=1}^{p_1+\cdots+p_m-k} t^q \sum_{\substack{l_1 \le l_2 \le \cdots \le l_q \\ l_1+\cdots+l_q = p_1+\cdots+p_m-k}} h_q(l_1,\dots,l_q) \prod_{s=1}^{q} C^{(l_s)}. \tag{16}$$

Here, $h_q(l_1,\dots,l_q)$ is a combinatorial factor given by the number of partitions of $l_1 + \cdots + l_q$ objects in subsets with $l_1, l_2, \dots, l_q$ elements. It can be explicitly calculated as follows. Let $\{l_s : s = 1,\dots,q\} = \{i_1,\dots,i_s\}$ be the set of different values of the $l_1,\dots,l_q$ and let $q_u = \sharp\{n = 1,\dots,q : l_n = i_u\}$, $u = 1,\dots,s$. Then

$$h_q(l_1,\dots,l_q) = \binom{l_1+\cdots+l_q}{l_1 \cdots l_q} \frac{1}{q_1! \cdots q_s!}. \tag{17}$$

The expressions (14) and (16) simplify further if certain coefficients $C^{(n)}$ and $\lambda^{(p)}$ vanish, for example, because of symmetries. As this is easier to understand in terms of graphs, we postpone the discussion of this point to the following subsection.

In the case that $f$ is a probability density and $\Psi$ is of the form (5), it is clear that $\Phi_t$, as a convolution of a probability measure with a probability density, is a probability density for all $t$. Note that the first $N$ terms of the expressions (14) or (16) give a polynomial in $\phi$ and $t$. Obviously, this plays havoc with normalization conditions. We thus have to keep in mind that the approximation given by the first $N$ terms in the perturbation series is only good for $t$ and $\phi$ sufficiently small. In Sections 2.4, 3.1 and 3.2, we will discuss strategies to deal with this problem.

On the other hand, the expansion in terms of polynomials in $\phi$ makes it simplest to integrate the probability density and to get a formula for the distribution function. Taking, for example, the one-dimensional formula (16), we obtain, in the sense of asymptotic series,

$$P(a < Z_t \le b) \\ = \sum_{m=0}^{\infty} \frac{(-1)^m}{m!} \sum_{p_1,\dots,p_m=0}^{\bar{p}} \prod_{l=1}^{m} \lambda^{(p_l)} \sum_{k=0}^{p_1+\cdots+p_m} \binom{p_1+\cdots+p_m}{k} \frac{a^{k+1} - b^{k+1}}{k+1} \tag{18}$$

$$\times \sum_{q=1}^{p_1+\cdots+p_m-k} t^q \sum_{\substack{l_1\le l_2\le\cdots\le l_q \\ l_1+\cdots+l_q=p_1+\cdots+p_m-k}} h_q(l_1,\ldots,l_q)\prod_{s=1}^{q} C^{(l_s)},$$

where $Z_t$ is a random variable with probability density function $\Phi_t$, that is, a Lévy process with generator $\Psi$ and initial distribution $f$. Again, the first $N$ terms of the right-hand side give a good approximation for $a, b, t$ and $\lambda^{(p)}$ sufficiently small.

## 2.3. The Feynman graph representation

In this subsection, we represent the solution of the Cauchy problem (14) as a sum over generalized Feynman graphs that are evaluated by a simple algorithm called "Feynman rules." In physics, the result is known as the Feynman graph representation for amputated Green's functions for the case when the generator is some kind of Laplacian, that is, $\Psi(\nabla) = \sum_{X_1,X_2=1}^{d} D_{X_1,X_2}\frac{\partial^2}{\partial\phi_{X_1}\partial\phi_{X_2}}$, [18]. Combining this with the generalized Feynman graph calculus developed in [5, 6], we obtain the first main result of this article.

A graph is a geometrical object which consists of vertices and edges, that is, lines that connect exactly two vertices. Here, we need graphs with non-directed edges that have vertices of different types. Vertices can be distinguishable or non-distinguishable, depending on their type, and they can have distinguishable or non-distinguishable legs. We use the term "leg" for the part of the edge meeting the vertex. All of these graph theoretic notions are properly defined in Appendix A.2. But for a first reading, an intuitive comprehension is sufficient.

Our expansion is represented as a sum over generalized Feynman graphs. They have the following properties.

***Definition 2.8.*** *A generalized amputated Feynman graph is a graph with three types of vertices, called* inner full,[2] inner empty *and outer empty vertices, respectively; see Table 1. By definition full vertices are distinguishable and have distinguishable legs, whereas empty vertices are non-distinguishable and have non-distinguishable legs.*[3] *Outer empty are met by one edge only. Edges are non-directed and connect full and empty (inner and outer) vertices, but never connect two full or two empty vertices.*

*Let $m \in \mathbb{N}_0$. The set of generalized amputated Feynman graphs with $m$ inner full vertices with $p_1,\ldots,p_m$ the number of legs of the inner full vertices such that $p_j \le \bar{p}$ and $\lambda^{(p_j)} \neq 0$, $j = 1,\ldots,m$, is denoted by $\bar{F}(m)$.*

The name *generalized amputated Feynman graphs* has been chosen to distinguish the graphs used here from classical Feynman graphs, but also from the generalized Feynman graphs with outer full vertices used in [5, 6] (this latter type gives the right combinatorics to calculate the moments of $\Phi_t$). To simplify the notions, we here use the term Feynman

[2]In physics this type is called an interaction vertex.

[3]Empty inner vertices can be distinguished from empty outer vertices, however.

graph when we mean a generalized amputated Feynman graph and there is no danger of confusion.

Let $m, p_1, \ldots, p_m \in \mathbb{N}_0$. We construct a one to one correspondence between pairs $(K; I)$, where $K \subseteq \Omega(p_1, \ldots, p_m)$ and $I \in \mathcal{P}(\Omega(p_1, \ldots, p_m) \setminus K)$, and the Feynman graphs $\bar{F}(m, p_1, \ldots, p_m)$ with $m$ inner full vertices with $p_1, \ldots, p_m$ legs in the following way:

- For $j = 1, \ldots, m$, label the legs of the $j$th inner full vertex $(1, j), \ldots, (p_j, j)$. In this way, we obviously obtain a one-to-one correspondence between the legs of the inner full vertices and the elements in $\Omega(p_1, \ldots, p_m)$.

We first consider the case where the Feynman graph $G$ is given and we construct $(K, I)$:

- Define $K \subseteq \Omega(p_1, \ldots, p_m)$ as the set corresponding (under the above bijection) to those legs of inner full vertices in the Feynman graph $G$ that run into an outer empty vertex.
- Let $k$ be the number of inner empty vertices in $G$. For $q = 1, \ldots, k$, let $I_q$ be the subset in $\Omega(p_1, \ldots, p_m)$ corresponding to those legs of inner full vertices that run into the $q$th inner empty vertex. Finally, let $I = \{I_1, \ldots, I_k\}$. Then, evidently, $I \in \mathcal{P}(\Omega(p_1, \ldots, p_m) \setminus K)$.

Conversely, let $K \subseteq \Omega(p_1, \ldots, p_m)$ and $I \in \mathcal{P}(\Omega(p_1, \ldots, p_m) \setminus K)$, $I = \{I_1, \ldots, I_k\}$, be given. We construct a Feynman graph $G \in \bar{F}(m, p_1, \ldots, p_m)$:

- Draw $m$ inner full vertices with $p_1, \ldots, p_m$ legs and $k$ inner empty vertices with $\sharp I_1, \ldots, \sharp I_k$ legs and $p = \sharp K$ outer empty vertices.
- Connect the legs of the $q$th inner empty vertex with the legs of the inner full vertices corresponding to the elements in $I_q$.
- Also connect the legs of outer empty vertices with the legs of the inner full vertices corresponding to the elements of $K$. The result is a Feynman graph $G \in \bar{F}(m)$.

For an illustration of the above construction, see Figure 1. We have established the following result.

**Lemma 2.9.** *Let $m, p_1, \ldots, p_m \in \mathbb{N}_0$. There then exists a one-to-one correspondence between $\{(K, I) : K \subseteq \Omega(p_1, \ldots, p_m), I \in \mathcal{P}(\Omega(p_1, \ldots, p_m) \setminus K)\}$ and $\bar{F}(m, p_1, \ldots, p_m)$, the set of Feynman graphs with $m$ inner full vertices with $p_1, \ldots, p_m$ legs.*

The following algorithm in physics is named "Feynman rules" and gives a numerical value to each Feynman graph.

**Table 1.** Different types of vertices

| | **Full** | **Empty** |
|---|---|---|
| Inner | $\bullet$ | $\circ$ |
| Outer | $\times$ | $\otimes$ |

***Definition 2.10.*** *Let $G \in \bar{F}(m)$ be such that $p_1, \ldots, p_m$ are the numbers of legs of the inner full vertices of $G$. Then, for $\phi \in \mathbb{R}^d$ and $t \in (0, \infty)$, $\mathcal{V}[G] = \mathcal{V}[G](t, \phi)$ is a real number obtained by the following method*[4]*:*

1. *to each leg of an inner full vertex, assign a value $X_a = X_j^l \in \{1, \ldots, d\}$, with $a$ the element in $\Omega(p_1, \ldots, p_m)$ corresponding to that leg;*
2. *multiply each outer empty vertex with a factor $(-\phi_{X_b})$, where $b \in \Omega(p_1, \ldots, p_m)$ is the element corresponding to a leg of the inner full vertex which is directly connected with that outer empty vertex;*
3. *for each inner empty vertex with $l$ legs, multiply with a factor $tC_J = tC_{X_{a_1} \cdots X_{a_l}}$, where $J = \{a_1, \ldots, a_l\} \subseteq \Omega(p_1, \ldots, p_m)$ is the subset corresponding to the $l$ legs of inner full vertices that are directly connected to the inner empty vertex;*
4. *for the $j$th inner full vertex, $j = 1, \ldots, m$, multiply with a factor $\lambda^{(p_j)}_{X_1^j \cdots X_{p_j}^j}$ and sum up over $\sum_{X_1^j, \ldots, X_{p_j}^j = 1}^d$.*

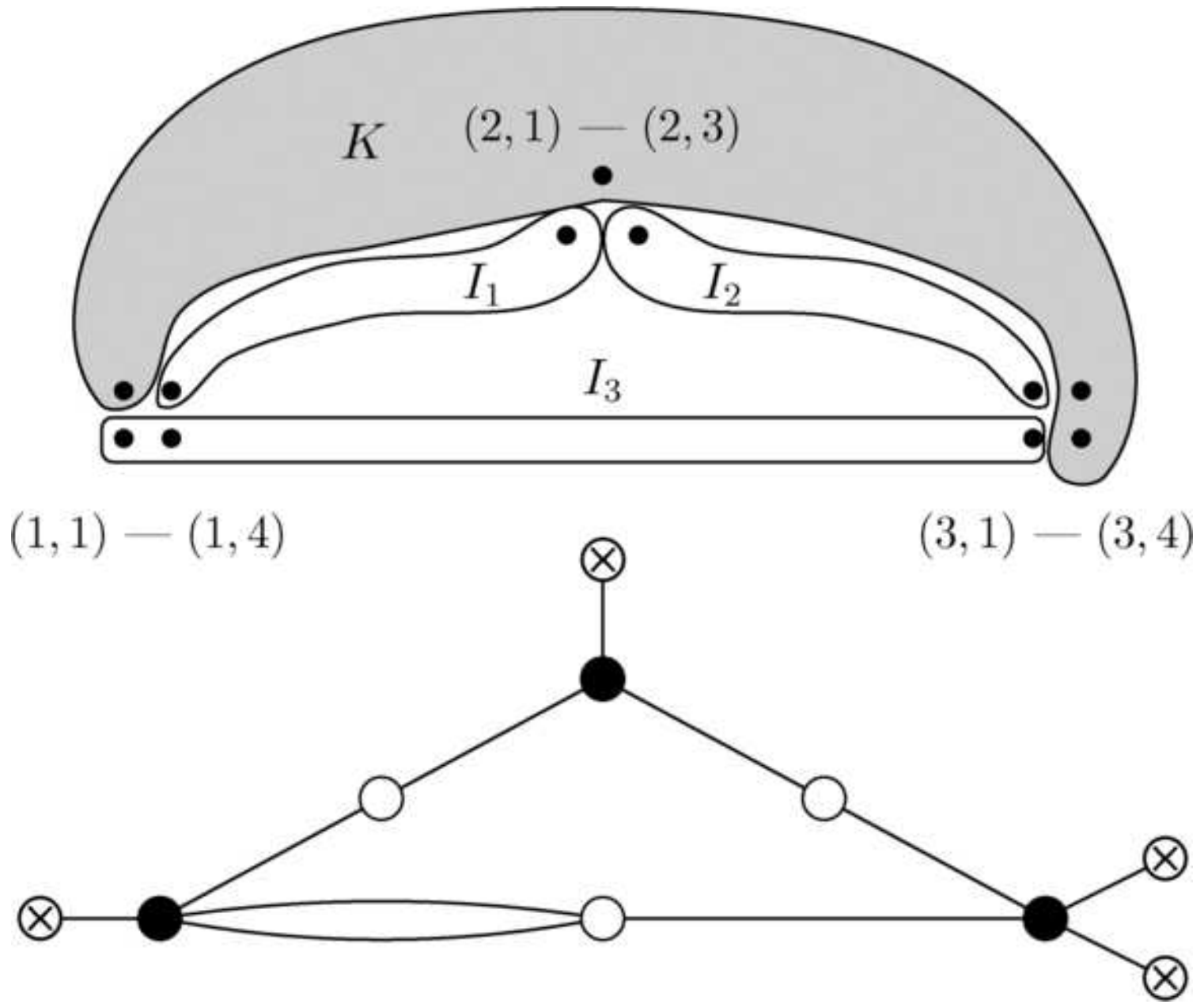

**Figure 1.** Construction of a generalized Feynman graph from the set $K$ and the partition $I = \{I_1, I_2, I_3\}$.

[4]We use the labeling of legs of inner full vertices by $\Omega(p_1, \ldots, p_m)$ introduced above, but any other labeling would do.

Combining Proposition 2.7, Lemma 2.9 and Definition 2.10, one thus gets the Feynman graph expansion of the convolution semigroup.

**Theorem 2.11.** *The solution of the Cauchy problem (1) has the following asymptotic expansion:*

$$\Phi_t(\phi) = \sum_{m=0}^{\infty} \frac{(-1)^m}{m!} \sum_{G \in \bar{F}(m)} \mathcal{V}[G](t,\phi). \tag{19}$$

Apart from the more efficient notation that comes with the use of Feynman graphs, Feynman graphs can be used to reduce the number of terms in the expansion. Note that the value $\mathcal{V}[G](t,\phi)$ attributed to the graph $G$ depends only on the topological Feynman graph. Hence, one can change the sum in (19) into a sum over topological Feynman graphs by multiplying $\mathcal{V}[G](t,\phi)$ with a weight factor or multiplicity, that is, with the number of Feynman graphs that give the same topological Feynman graph. Such multiplicities in low orders can be calculated by hand; see (4) for an example.

Equation (19) is still fairly general. Often, a number of Feynman graphs do not give contributions. For example, in the case where $\Psi(\xi) = \Psi(-\xi)$, one has $C^{(n)} = 0$ for $n$ odd and hence all Feynman graphs with at least one empty vertex with an odd number of legs can be omitted in the expansion.

Another simplification occurs if the $V(\varphi)$ is quadratic in $\varphi$. Then, the only kind of inner full vertex that occurs is of the type —•— and can thus be considered as an edge of new type. Note that the legs of the inner full vertex are distinguishable, hence this new type of edge is directed, which changes multiplicities. Finally, if $\lambda^{(2)}_{X_1X_2} = \lambda\delta_{X_1X_2}$ is a multiple of the Kronecker symbol, the evaluation of a graph such as the one on the right-hand side of equation (4) can be simplified as follows.

**Corollary 2.12.** *If $V(\phi)$ is quadratic, then the set of Feynman graphs $\bar{F}(m)$ of order $m$ can be identified with the set of all graphs $\bar{Q}(m)$ with two types of indistinguishable vertices, called* outer empty *and* inner empty*, with $m$ directed edges such that an outer empty vertex is hit by exactly one edge.*

*If, in addition, $\lambda^{(2)}$ is diagonal,*[5] *the value $\mathcal{V}[G](t,\phi)$ attributed to such a graph $G \in \bar{Q}(m)$ can be calculated according to the following rules:*

1. *to each of the $m$ edges, assign a value $X_a \in \{1,\ldots,d\}$, $a = 1,\ldots,m$;*
2. *for each inner empty vertex with $l$ legs, multiply with a factor $tC_{X_{a_1}\cdots X_{a_l}}$, where $a_1,\ldots,a_l$ correspond to the numbers given to the edges that are connected to this vertex*[6]*;*
3. *for each edge, add a factor $\lambda$, or, alternatively, add an overall factor $\lambda^m$;*
4. *sum up over all edges, that is, perform the sum*[7] $\sum_{X_1,\ldots,X_m=1}^{d}$.

[5] If $\lambda^{(2)}$ is symmetric, this can always be achieved by choosing a suitable basis on $\mathbb{R}^d$.

[6] If an edge connects the vertex with itself, the corresponding index has to be repeated twice.

[7] Of course, for $d = 1$, the same simplifications as discussed at the end of Section 2.2 occur.

As a side remark, we note that the combinatorics of the Feynman graph formalism are dimension independent. Only Feynman rules depend on dimension. The number of steps needed in the evaluation of a given Feynman graph grow like $\sim d^m$ with the dimension, where $m$ is the order of the graph. If one compares this with an evaluation of $\Phi_t$ by a multidimensional Fourier transform, the steps needed for the evaluation grow rather like $\sim m^d$, with $m$ a parameter presenting, for example, the number of summands in an approximation of the one-dimensional Fourier integrals by a sum. For large dimensions, the Feynman graph formalism thus promises better results than numerical evaluation of Fourier transforms. Since the Feynman graph formalism has its roots in infinite dimensions, this should not be a surprise.

## 2.4. The logarithm of the solution

In this subsection, we obtain an asymptotic expansion of the logarithm of the solution in terms of connected Feynman graphs. In many points, the argument follows [5].

A graph $G \in \bar{F}(m)$ is said to be *connected* if and only if each two vertices are *connected* by a walk passing through finitely many edges. We denote the collection of connected $m$th order Feynman graphs by $\bar{F}_c(m)$.

On the level of partitions, connectedness can be expressed as follows.

***Definition 2.13.*** *Let $m, p_1, \dots, p_m \in \mathbb{N}$ be given and $J_j = \{(j,1), \dots, (j,p_j)\}$, $j = 1, \dots, m$. We say that the pair $(K, I)$ with $K \subseteq \Omega(p_1, \dots, p_m)$, $I \in \mathcal{P}(\Omega(p_1, \dots, p_m) \setminus K)$, $I = \{I_1, \dots, I_k\}$, is connected if it fulfills the following condition: $\nexists 1 \leq i_1, \dots, i_q \leq k$, $1 \leq q < k$ and $1 \leq j_1, \dots, j_s \leq m$, $1 \leq s < m$, such that*

$$\bigcup_{\alpha=1}^{q} I_{i_\alpha} = \bigcup_{\alpha=1}^{s} (J_{j_\alpha} \setminus K). \tag{20}$$

Let $\mathcal{P}_c(\Omega(p_1, \dots, p_m) \setminus K)$ be the collection of partitions $I \in \mathcal{P}(\Omega(p_1, \dots, p_m) \setminus K)$ such that $(K, I)$ is connected.

**Lemma 2.14.** *A graph $G \in \bar{F}(m)$ is connected if and only if the corresponding pair $(K, I)$ is connected.*

**Proof.** We only sketch the idea of the proof. In terms of the corresponding Feynman graph $G$, the condition in Definition 2.13 means that one cannot find a true subset $Q_1, Q_2$ of empty and full vertices such that all edges from an empty vertex in $Q_1$ end up at a full vertex in $Q_2$ and vice versa. This is evidently equivalent to the connectedness of $G$. For details, see [5]. □

***Definition 2.15.*** *Let $J_j$, $j = 1, \dots, m$, be as in Definition 2.13, $\phi \in \mathbb{R}^d$ and $X_a = X_j^l \in \{1, \dots, d\}$ for $a = (j,l) \in \Omega(p_1, \dots, p_m)$. We use the abbreviation $J_j^\phi = \prod_{a \in J_j} (\varphi_{X_a} - \phi_{X_a})$,* $j = 1,\dots,m$. *The block truncated moment functions* $\langle\prod_{l=1}^m J_l^\phi\rangle_{\nu_t}^{(T)}$ *of* $\nu_t$ *are recursively (in* $m\in\mathbb{N}$*) defined as follows:*

$$\langle J_1^\phi\cdots J_m^\phi\rangle_{\nu_t} = \sum_{\substack{I\in\mathcal{P}(1,\dots,m)\\ I=\{I_1,\dots,I_k\}}}\prod_{l=1}^k\left\langle\prod_{q\in I_l}J_q^\phi\right\rangle_{\nu_t}^{(T)}. \tag{21}$$

The symbol $(T)$ in Definition 2.15 means that each polynomial (in $\varphi$) $J_q^\phi$ in the combinatorics of truncation is treated as a single object.

**Proposition 2.16.** *Let* $J_1,\dots,J_m$ *be as in Definition 2.13. We then obtain the following expansion of block truncated moments into ordinary truncated moments:*

$$\langle J_1^\phi\cdots J_m^\phi\rangle_{\nu_t}^{(T)} = \sum_{\substack{K\subseteq\Omega(p_1,\dots,p_m)\\ I\in\mathcal{P}_c(\Omega(p_1,\dots,p_m)\setminus K)}}\prod_{j\in K}(-\phi_{X_j})\prod_{l=1}^k\langle I_l\rangle_{\nu_t}^T. \tag{22}$$

**Proof.** To verify equation (22), we have to insert the right-hand side of this equation into the defining equation for the left-hand side. We thus have to show that

$$\langle J_1^\phi\cdots J_m^\phi\rangle_{\nu_t} = \sum_{\substack{I\in\mathcal{P}\{1,\dots,m\}\\ I=\{I_1,\dots,I_k\}}}\prod_{l=1}^k\left[\sum_{\substack{K_l\subseteq\bigcup_{q\in I_l}J_q\\ Q_l\in\mathcal{P}_c(\bigcup_{q\in I_l}J_q)\\ Q_l=\{Q_{l,1},\dots,Q_{l,k_l}\}}}\prod_{j\in K_l}(-\phi_{X_j})\prod_{s_l=1}^{k_l}\langle Q_{l,s_l}\rangle_{\nu_t}^T\right]. \tag{23}$$

On the other hand, $\langle J_1^\phi\cdots J_m^\phi\rangle_{\nu_t}$ can be directly expanded in terms of truncated moments:

$$\langle J_1^\phi\cdots J_m^\phi\rangle_{\nu_t} = \sum_{\substack{K\subseteq\Omega(p_1,\dots,p_m)\\ R\in\mathcal{P}(\Omega(p_1,\dots,p_m)\setminus K)\\ R=\{R_1,\dots,R_k\}}}\prod_{j\in K}(-\phi_{X_j})\prod_{l=1}^k\langle R_l\rangle_{\nu_t}^T. \tag{24}$$

We have to prove the equality between the right-hand side of (23) and (24), that is, we have to construct a correspondence between the pairs $(K,R)$ and the objects $(I,(K_1,Q_1),\dots,(K_k,Q_k))$ indexing the sum in (23) such that their contributions to the sum are equal.

Let $(R,K)$ be given. We say that $R$ *connects* $q$ and $j$ (in notation, $q\sim_R j$) if the full inner vertices corresponding to $J_q$ and $J_j$, respectively, are connected in the generalized Feynman graph corresponding to $R$. It is easy to see that $R$ is an equivalence relation on $\{1,\dots,m\}$. Let $I=\{I_1,\dots,I_k\}$ be the equivalence class of $\sim_R$. Then $I\in\mathcal{P}(\{1,\dots,m\})$. We set $K_l=(\bigcup_{q\in I_l}J_q\cap K)$, $l=1,\dots,k$, and $Q_l=\{S\in R: S\subseteq\bigcup_{q\in I_l}J_q\}$. Using the equivalence relation $\sim_R$, it is now easy to show that $Q_l \in \mathcal{P}_c(\bigcup_{q\in I_l} J_q)$; see [5] for the details.

The converse construction is trivial: take $K = \bigcup_{l=1}^k K_l$ and $R = \bigcup_{l=1}^k Q_l$. These two constructions give mappings between the index sets of the sums in (23) and (24) such that the summands are equal. Furthermore, the mappings are inverses of each other, which establishes the one-to-one correspondence. □

We have now completed the preparations to prove the following theorem.

**Theorem 2.17.** *In the sense of asymptotic series,* $\log \Phi_t(\phi)$ *has the following expansion in terms of connected Feynman graphs:*

$$\log \Phi_t(\phi) = \sum_{m=1}^{\infty} \frac{(-1)^m}{m!} \sum_{G\in \bar{F}_c(m)} \mathcal{V}[G](t,\phi). \tag{25}$$

**Proof.** Note that by the ordinary linked cluster theorem (Theorem 2.6), the following holds, in the sense of asymptotic series:

$$\log \Phi_t(\phi) = \log\langle \nu_t, \mathrm{e}^{-V_\phi}\rangle = \sum_{m=1}^{\infty} \frac{(-1)^m}{m!} \langle V_\phi^m \rangle_{\nu_t}^{(T)}. \tag{26}$$

The assertion thus follows by application of Proposition 2.16 and Lemma 2.14. □

The above result is useful in three respects.

First, if $\Psi$ is the symbol of a Lévy process and the initial condition $f = \mathrm{e}^{-V}$ is a probability density, then the solution $\Phi_t(\phi)$ is a probability density for all $t$. We have seen in Proposition 2.7 and in equation (16) that the first $N$ summands of the asymptotic expansion for $\Phi_t(\phi)$ is a polynomial in $\phi$ and thus non-normalizable. However, we derived an asymptotic expansion for $\log \Phi_t(\phi)$ such that the first $N$ terms are a polynomial $P_N(\phi)$ in $\phi$ that goes to $-\infty$ for large $\phi$ (this is often the case for $N$ odd) and we obtain $\Phi_t(\phi) \approx \mathrm{e}^{-P_N(\phi)}$ as the solution in $N$th order perturbation theory. While the large $\phi$ behavior of this approximate solution is obviously still problematic,[8] one has at least gained normalizability.

Second, in statistical estimations by the maximum likelihood method, one needs an explicit formula for $\log \Phi_t(\phi)$, which is just the sum over connected Feynman graphs.

Finally, the connectedness condition also reduces the number of terms in the expansion. For example, equation (4), only the 2nd, 3rd, 4th and 7th terms contribute to the logarithm of the solution.

[8]The probability that a random variable $Z_t$ distributed according $\Phi_t$ takes large values $\phi$ is suppressed much too strongly. To adjust the large $\phi$ behavior, this first step has to be combined with resummation techniques (cf. Section 3.2).

## 2.5. Borel summability

In certain cases, a function may be recovered from its asymptotic series expansion. This class of functions comprises the Borel summable functions. A sufficient condition for a function to equal its Borel sum is the following version of Watson's theorem due to Sokal [20].

**Theorem 2.18.** *Let $h$ be analytic in the interior of the circle $C_R := \{\zeta \in \mathbb{C} : \Re\zeta^{-1} > R^{-1}\}$ (cf. Figure 2(a)) such that*

$$\left| h(\zeta) - \sum_{m=0}^{N-1} a_m \zeta^m \right| \leq A\sigma^N N! |\zeta|^N \tag{27}$$

*uniformly in $N$ and $\zeta \in C_R$, with constants $A, \sigma > 0$. Under this assumption, one may show that $B(\tau) := \sum_{m=0}^{\infty} a_m \tau^m / m!$ converges for $|\tau| < 1/\sigma$ and has an analytic continuation to the striplike region $S_\sigma = \{\tau \in \mathbb{C} : \mathrm{dist}(\tau, \mathbb{R}_+) < 1/\sigma\}$ (cf. Figure 2(b)) obeying the bound*

$$|B(\tau)| \leq K \exp(|\tau|/R) \tag{28}$$

*uniformly in every $S_{\sigma'}$, with $\sigma' > \sigma$. Furthermore, $h$ can be represented by the absolutely convergent integral*

$$h(\zeta) = \frac{1}{\zeta} \int_0^\infty \mathrm{e}^{-\tau/\zeta} B(\tau)\, \mathrm{d}\tau. \tag{29}$$

*Conversely, if $B(\tau)$ is a function analytic in $S_\sigma$ and there satisfying (28), then the function defined by (29) is analytic in $C_R$ and satisfies (27) with $a_m = \frac{\mathrm{d}^m}{\mathrm{d}\tau^m} B(\tau)|_{\tau=0}$ uniformly in every $C_{R'}$ with $R' < R$.*

If $h$ is analytic on $C_\infty = \bigcup_{R>0} C_R = \{\zeta \in \mathbb{C} : \Re\zeta > 0\}$ and (27) holds on that domain, then of course the representation (29) also holds for $\zeta \in C_\infty$ and for $\zeta$ real and $\beta = \Re\zeta \in \mathbb{R}_+$ in particular. The aim now is to apply this result to the series given in (19). It is easy to show analyticity, as seen in the following result.

**Lemma 2.19.** *The function $\mathbb{R}_+ \ni \beta \to \Phi_t^\beta(\phi)$ (cf. Section 2.1) has analytic continuation $C_\infty \ni \zeta \to \Phi_t^\zeta(\phi)$, $\forall t \in \mathbb{R}_+, \phi \in \mathbb{R}^d$.*

**Proof.** Note that for $\zeta \in C_\infty$, $f_\zeta = \mathrm{e}^{-\zeta V} \in \mathcal{S}(\mathbb{R}^d)$ and hence $\Phi_t^\zeta(\phi) = \nu_t * f_\zeta(\phi)$ is well defined. To show complex differentiability in $\zeta$, we now proceed as in the proof of Lemma A.1. □

In order to check (27), we need an estimate on the multiple $\zeta$ derivatives of $\Phi_t^\zeta(\phi)$ on $C_R$ or even on $C_\infty$. As $\Phi_t^\zeta(\phi)$ is not known explicitly, such an estimate is difficult to obtain. But, in the case that $\Psi$ is the symbol of a Lévy process, the probabilistic structure helps.

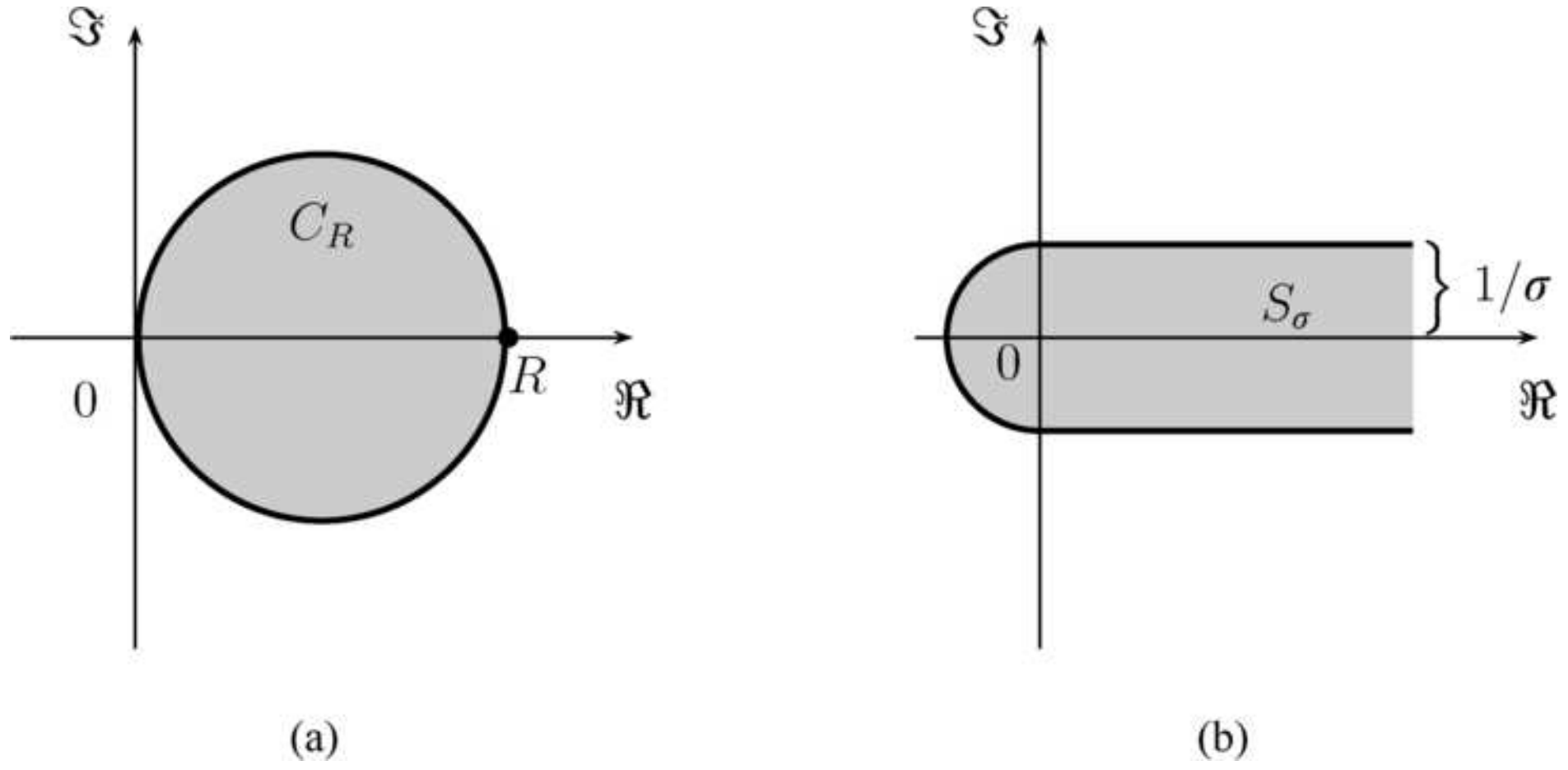

**Figure 2.** (a) The ball $C_R$ and (b) the striplike region $S_\sigma$.

**Lemma 2.20.** *If* $\Psi$ *is the symbol of a Lévy process and* $V(\varphi) \geq 0\ \forall \varphi \in \mathbb{R}^d$, *then*

$$\left| \frac{\mathrm{d}^m}{\mathrm{d}\zeta^m} \Phi_t^\zeta(\phi) \right| \leq \langle V_\phi^m \rangle_{\nu_t}, \qquad \zeta \in C_\infty, t \in \mathbb{R}_+, \phi \in \mathbb{R}^d. \tag{30}$$

**Proof.** We have $|\frac{\mathrm{d}^m}{\mathrm{d}\zeta^m} \Phi_t^\zeta(\phi)| = |\langle \nu_t, V_\phi^m \cdot \mathrm{e}^{-\zeta V_\phi} \rangle| \leq \langle V_\phi^m \rangle_{\nu_t}$ as $\nu_t$ is a probability measure and $|\mathrm{e}^{-\zeta V_\phi}| \leq 1$. □

Lemma 2.20 allows us, in the given situation, to check (27) by Taylor's estimate for the remainder and the calculation of the $N$th coefficient in the perturbation series. Being a sum over Feynman graphs, the coefficient can be dominated by the number of such Feynman graphs times the maximal value of $\mathcal{V}[G](t, \phi)$ over all these graphs. The latter is easy to control. We therefore start with an estimate on the number of Feynman graphs. Let us first assume that the shift $\phi$ is zero. We furthermore assume that $V$ is homogeneous of degree $\bar{p}$. Then, only graphs without outer empty vertices and with inner full vertices with $\bar{p}$ legs contribute. As such graphs are in one-to-one correspondence with the number of partitions of $\bar{p}N$ objects, we have to control the number of such partitions.

The number of partitions of $m$ objects is known as $m$th *Bell number* and we denote it by $b_m$. It has the asymptotic behavior (see, e.g., [12])

$$b_m \sim m^{-1/2} \lambda(m)^{m+1/2} \mathrm{e}^{\lambda(m)-m-1}, \tag{31}$$

where $\lambda(m)$ is implicitly defined by $\lambda(m) \ln \lambda(m) = m$. Using Stirling's formula $N! \sim \sqrt{2\pi} \mathrm{e}^{(N+1/2)\ln N - N}$, we get the following asymptotic formula for $b_{\bar{p}N}/N!$:

$$\frac{1}{N!} b_{\bar{p}N} \sim \frac{1}{\sqrt{2\pi}} \bar{p}^{-1/2} \mathrm{e}^{(\bar{p}-1)(N \ln N - N) - 1/2 \ln N - 1} \mathrm{e}^{\bar{p}N(\ln \bar{p} - \ln \ln \lambda(\bar{p}N))} \mathrm{e}^{(1/2)\ln \lambda(\bar{p}N) + \lambda(\bar{p}N)}, \tag{32}$$

from which we obtain that

$$\frac{1}{N!} b_{\bar{p}N} \leq A' K^N N! \tag{33}$$

is fulfilled for $\bar{p} \leq 2$ and $K$ sufficiently large. In the following, we restrict ourselves to $\bar{p} = 2$.

Let us now consider a Feynman graph $G \in \bar{F}(N)$ without empty outer vertices. If the modulus of the coefficients $C^{(n)}_{X_1,\ldots,X_n}$ is uniformly bounded by $\hat{A}\kappa^n$ for $\hat{A}$, $\kappa$ sufficiently large, we obtain

$$|\mathcal{V}[G](t,0)| \leq \max\{\hat{A},1\}^N \hat{\kappa}^N d^N \kappa^{2N} \max\{1,t\}^{2N} \qquad \forall G \in \bar{F}(N), \tag{34}$$

with $\hat{\kappa} = \max_{X_1,X_2 \in \{1,\ldots,d\}} |\lambda^{(2)}_{X_1 X_2}|$. Thus, under the given conditions, the estimate (27) is fulfilled for $\phi = 0$ with $h(\zeta) = \Phi^{\zeta}_t(0)$, $a_m = \frac{(-1)^m}{m!} \sum_{G \in \bar{F}(m)} \mathcal{V}[G](t,0)$, $A = A'$ and $\sigma = d\hat{\kappa}\kappa^2 K \max\{1,t\}^2 \max\{\hat{A},1\}$. It turns out that the case $\phi \neq 0$ can be easily traced back to the case $\phi = 0$.

**Theorem 2.21.** *Let $\Psi$ be a Lévy symbol as in (5) such that* $\operatorname{supp}(r)$ *is compact and let $f = \mathrm{e}^{-V}$ with $V(\varphi) = \langle \lambda^{(2)}, \varphi^{\otimes 2} \rangle_2$. Then $\Phi^{\zeta}_t(\phi)$, as a function of $\zeta$, fulfills the conditions of Theorem 2.18. In particular, the series*

$$B(\tau,t,\phi) = \sum_{m=0}^{\infty} \frac{(-\tau)^m}{(m!)^2} \sum_{G \in \bar{F}(m)} \mathcal{V}[G](t,\phi) \tag{35}$$

*converges on a ball of radius $1/\sigma$ centered at 0 and has analytic continuation to the strip $S_\sigma$ (for some $\sigma = \sigma(t,\phi) < \infty$). Denoting this continuation by the same symbol, we get the solution $\Phi_t$ of the Cauchy problem (1) via the formula*

$$\Phi_t(\phi) = \int_0^\infty \mathrm{e}^{-\tau} B(\tau,t,\phi)\, \mathrm{d}\tau. \tag{36}$$

**Proof.** If $\operatorname{supp}(r) \subseteq \{x \in \mathbb{R}^d : |x| \leq \kappa\}$, then the uniform bound $|C^{(n)}_{X_1 \cdots X_n}| \leq \hat{A}\kappa^n$ holds for $\hat{A}$ sufficiently large, since $C^{(n)}_{X_1 \cdots X_n}$ for $n > 2$ up to a factor $z$ is an $n$th moment of $r$. From the considerations above, the assertion now follows for $\phi = 0$.

If $\phi \neq 0$, we modify for a fixed time $t = t_0$ the Lévy symbol $\Psi$ by a deterministic term $\Gamma_{\phi,t_0}(-\mathrm{i}\xi) = \Psi(-\mathrm{i}\xi) + \mathrm{i}\phi \cdot \xi / t_0$ generating the shift $\phi$ at the time $t_0$. If we let $\Theta_{t,t_0,\phi}$ be the solution of the Cauchy problem with $\Psi$ replaced by $\Gamma_{t_0,\phi}$, then obviously $\Phi^{\zeta}_{t_0}(\phi) = \Theta^{\zeta}_{t_0,t_0,\phi}(0)$. Since the coefficients of $\Psi$ fulfill the above uniform bound, the same holds for the coefficients of $\Gamma_{t_0,\phi}$ (with $\hat{A}$ changed depending on $t_0$ and $\phi$). The assertion holds for $\Theta^{\zeta}_{t_0,t_0,\phi}(0)$, hence it also holds for $\Phi_{t_0}(\phi)$.

In fact, the two expansions are equal. The only difference is that in the graphic expansion of $\Theta_{t_0,t_0,\phi}(0)$, there are no empty outer vertices. This is compensated for by evaluation of the empty inner vertices with one leg by $(C^{(1)}_X - \phi_X)$. Multiplying out the parenthesis just generates the contributions of a non-modified inner empty vertex with one leg and the outer empty vertex. □

The restriction to $\bar{p} = 2$ still covers a set of initial conditions $f$ that determine the convolution kernel $\nu_t$ completely (cf. Section 2.1). It is also sufficient for the most important application that we have in mind; see the following Section 3.1.

The second technical restriction on the support of the jump distribution $r$ is more severe. It can be overcome by considering generalized notions of Borel summability (cf. [13]) that for the case of Lévy processes can be adapted to a growth like $(\sim N!)^\gamma$ of the moments of the jump distribution $r$. In applications, this matters if one wants to choose (36) or generalizations thereof as the point of departure of the numerical evaluation (cf. [10], Chapters 16.3–16.6) and one has to take care concerning the speed of convergence of the (generalized) Borel transform. This program, as with all high precision numerics based on our proposal, is beyond the scope of the present article.

It would be desirable to have a similar theorem for $\log \Phi_t(\phi)$. The analyticity of $\log \Phi_t^\zeta(\phi)$ on some circle $C_R$ is immediate. The bounds (27) for the perturbation series *at* zero is a corollary of the above considerations since the connectedness condition reduces the total number of Feynman graphs. In particular, the Borel sum $B_c(\tau, t, \phi)$ of $\log \Phi_t^\zeta(\phi)$ has a non-vanishing radius of convergence. But, for the log-solution, the analog of Lemma 2.20 is missing and therefore the bounds on the series coefficients *at* $\zeta = 0$ do not produce a uniform estimate (27) on some $C_R$, $R > 0$. So, the existence of an analytic continuation of $B_c(\tau, t, \phi)$ is not proved. Of course, one can repeatedly use the chain rule to obtain a bound on $\log \Phi_t^\zeta(\phi)$, but that generates an additional factor $\sim b_N$. So, more sophisticated estimates or more powerful versions of Theorem 2.18 are needed [13].

Here, we do not enter into this subject. Numerical procedures to evaluate the integral transform (29) are based on the convergence of the Borel sum which permits resummation as, for example, the conformal mapping method [10], Chapter 16. Such methods, relying on the first $N$ terms of the perturbation expansion only, automatically produce analytic continuations of an approximated $B_c(\tau, t, \phi)$ on some striplike region that are sufficiently bounded in order to render (29) absolutely convergent.

# 3. Application to Lévy distributions

## 3.1. Large diffusion expansion of Lévy densities

It is more interesting to get explicit (approximate) formulae for the distribution of a Lévy process starting at a fixed point (e.g., 0) than to determine the density of the distribution for an initial condition $f = \mathrm{e}^{-V}$ with $V$ as above. Though the densities of Lévy distributions can be expressed in terms of more or less elementary functions for many interesting examples (see, e.g., [1], Chapter 1), or can at least be expanded in a series in other cases [17], a generic formula is missing.

In this section, we present asymptotic expansions of (non-normalized) Lévy densities and their logarithms, respectively, for the case where the generator of the process has a large diffusion constant $\mu > 0$. The formulae are generic in the sense that no detailed assumptions on the distribution of jumps, $r$, are required. That is, we consider the generator (6) with $D_{X_1X_2} = \mu\delta_{X_1X_2}$, $\delta_{X_1,X_2}$ being the Kronecker symbol, and we assume that $\mu$ is large and hence $1/\mu$ is a small parameter in which we would like to expand. As the initial value problem for the density of the jump-diffusion process starting at zero, we get, for $t, \phi \in (0,\infty) \times \mathbb{R}^d$,

$$\begin{aligned} \frac{\partial \Phi_t}{\partial t}(\phi) &= \mu\Delta\Phi_t(\phi) - \sum_{X=1}^{d} a_X \frac{\partial}{\partial \phi_X}\Phi_t(\phi) + z\int_{\mathbb{R}^d\setminus\{0\}} [\Phi_t(\phi+\varphi) - \Phi_t(\phi)]\,\mathrm{d}r(\varphi), \\ \Phi_0(\phi) &= \delta_0(\phi), \end{aligned} \tag{37}$$

with $\delta_0$ the Dirac measure with weight one at zero and $\Delta$ the Laplacian on $\mathbb{R}^d$. The solution, $\Phi_t^{\text{diff.}}$, of the diffusive part alone is, of course,

$$\Phi_t^{\text{diff.}}(\phi) = (4\pi\mu t)^{-d/2}\mathrm{e}^{-\|\phi\|^2/(4\mu t)} = \left(\frac{\beta}{\pi}\right)^{d/2}\mathrm{e}^{-\beta\|\phi\|^2}, \qquad \beta = 1/(4\mu t), t > 0. \tag{38}$$

Let $\nu_t = \nu_t^{\text{jump}}$ be the inverse Fourier transform of $\mathrm{e}^{t\Psi^{\text{jump}}(-\mathrm{i}\xi)}$, where $\Psi^{\text{jump}}$ is obtained from $\Psi$ letting $\mu = 0$. Obviously, we obtain the solution of (37) by $\Phi_t(\phi) = \nu_t * \Phi_t^{\text{diff.}}$. With $V(\phi) = \|\phi\|^2$, $f_\beta = \mathrm{e}^{-\beta V}$ one thus gets $(\frac{\beta}{\pi})^{-d/2}\Phi_t(\phi) = \nu_t * f_\beta(\phi)$. We note that the expansion of the right-hand side in powers of $1/\mu$ is equivalent to the expansion in powers of $\beta$ as the $m$th coefficient only differs by a factor $1/(4t)^m$. Furthermore, from now on, we omit the multiplication with factors $\beta$ in the evaluation rules for the Feynman graphs in order to obtain a more explicit $\beta$ dependence. Combining this with Corollary 2.12 and Theorem 2.17, we obtain the $1/\mu$ (large diffusion) expansion for the non-normalized (log-) density function.

**Theorem 3.1.** *In the sense of asymptotic series in $\beta$, we obtain the following large diffusion expansion of the non-normalized solution $(\frac{\beta}{\pi})^{-d/2}\Phi_t$, where $\Phi_t$ solves (37):*

$$\left(\frac{\beta}{\pi}\right)^{-d/2}\Phi_t(\phi) = \sum_{m=0}^{\infty}\frac{(-\beta)^m}{m!}\sum_{G\in\bar{Q}(m)}\mathcal{V}[G](t,\phi). \tag{39}$$

*Furthermore, again in the sense of asymptotic series, we have*

$$\log\left[\left(\frac{\beta}{\pi}\right)^{-d/2}\Phi_t(\phi)\right] = \sum_{m=1}^{\infty}\frac{(-\beta)^m}{m!}\sum_{G\in\bar{Q}_c(m)}\mathcal{V}[G](t,\phi), \tag{40}$$

*with $\bar{Q}_c(m)$ the set of connected graphs in $\bar{Q}(m)$.*

Combining of (39) with Theorem 2.21 yields the following convergent representation as a Borel transform.

**Corollary 3.2.** *Suppose that the support of the jump distribution $r$ is compact. Let $B(\tau, t, \phi)$ be the analytic continuation of $\sum_{m=0}^{\infty} \frac{(-\tau)^m}{(m!)^2} \sum_{G\in\bar{Q}(m)} \mathcal{V}[G](t,\phi)$ to some strip of the form $S_\sigma$. Then,*

$$\Phi_t(\phi) = \pi^{-d/2}\beta^{d/2-1}\int_0^\infty \mathrm{e}^{-\tau/\beta} B(\tau, t, \phi)\, \mathrm{d}\tau. \tag{41}$$

A comparison between (39)–(41) and other possible expansions known from statistical mechanics seems to be of interest. In particular, readers familiar with the theory of classical gases might find the convergent small $z$ ("low activity") expansion [16] more suited than the small $\beta$ ("high temperature") expansion. For a small $z$ expansion, a graphical representation by Meyer graphs also exists. But there is a decisive disadvantage from the point of view of applications. The small $\beta$ expansion produces expressions that directly depend on the quantities of statistical importance – $z$ and the moments $r_n$ of the jump distribution $r$. Hence, these quantities can be obtained from an empirical distribution by fits. This can be done to some extent nonparametrically, that is, without any a prior assumption on the jump distribution, as any sequence of moments of the jump measure depending on arbitrary parameter sets can be inserted into our formula. In the Meyer series, however, the jump distribution $r$ enters via non-trivial $r$ integrals over products of Meyer functions $f(s, s', \phi, \phi') = \mathrm{e}^{-ss'|\phi-\phi'|^2} - 1$ which, except for a few cases with $r$ particularly simple, can only be done numerically.[9]

Similar difficulties arise if one wants to solve the inverse Fourier integral of $\mathrm{e}^{t\Psi(-\mathrm{i}\xi)}$. Note that this integral is of oscillatory type, so even its numerics is not really trivial.

## 3.2. Second-order perturbation theory

In this section, we give an illustration of how to actually calculate the $1/\mu$ expansion of a Lévy density and how to improve the result using resummation techniques. We take (40) as our starting point. Calculations are carried through up to second order in $\beta$ only and we use the simplest possible resummation algorithm. Thus, the idea of this section is by no means to provide high-accuracy numerics, but rather to illustrate techniques up to second order of perturbation theory, that in higher orders will provide reasonable numerics. On the other hand, this does not mean that the second order calculation is not at all interesting from a numerical point of view. A more detailed discussion can be found below. Calculations to higher orders and more sophisticated resummation (see, e.g., [7, 9, 10]) are left for the future.

In this subsection, we assume, without loss of generality, that $\Psi^{\mathrm{jump}}(-\mathrm{i}\xi)$ has a vanishing first coefficient $C^{(1)}$. If this is not the case, then this can be achieved by a shift $\phi \to \phi - a$ (cf. the proof of Theorem 2.21).

[9]Note that the small-time expansion in [2, 17] is equivalent to the small $z$ expansion.

Under these conditions, the second-order perturbation expansion in $\beta$ for the left-hand side of (39) is given in (4). Extracting the connected graphs, we obtain

$$\log\left[\left(\frac{\beta}{\pi}\right)^{-d/2}\Phi_t(\phi)\right]$$

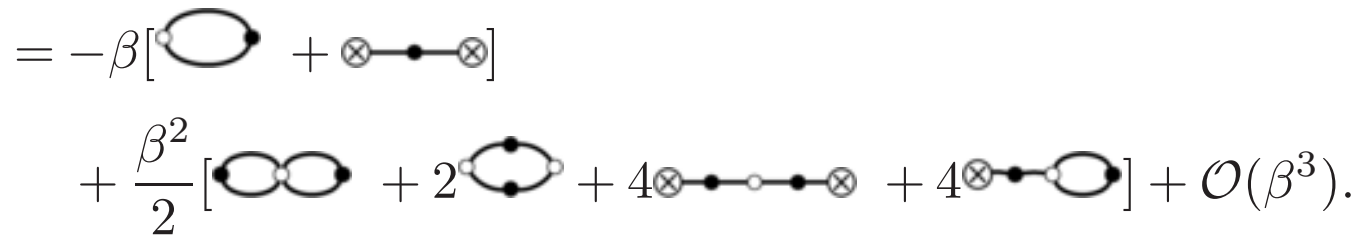

(42)

As in equation (4), we omitted the evaluation symbol $\mathcal{V}$ and multiplied by the multiplicities of the topological graphs.

Equation (42) should work for $z$ and $1/\mu$ sufficiently small and[10] $t \approx 1$, but the positive sign in front of the $\beta^2$ term implies that the exponential of the second-order approximation becomes non-normalizable (not $L^1$-integrable) if $\beta$ (or $z$) passes a certain threshold. To extract reliable data for higher values of $\beta$, resummation is required.[11] The idea of Padé approximants is to replace the Taylor polynomial by a rational function that has the same Taylor series as the polynomial. For definitions and details, we refer to Appendix A.3.

Following the standard calculation in A.3, we obtain as the second-order Padé approximant

$$\log\left[\left(\frac{\beta}{\pi}\right)^{-d/2}\Phi_t(\phi)\right] =$$

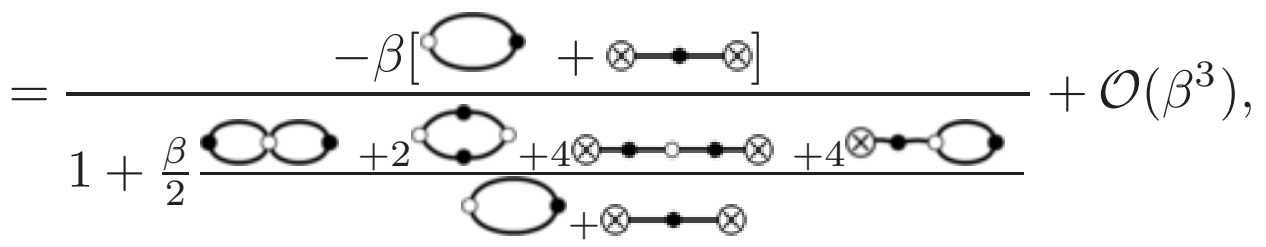

(43)

which, if exponentiated, gives a normalizable approximation to $\Phi_t$ for $t, z$ and $\beta$ arbitrary. Evaluation for $d = 1$ with $C^{(n)} = zr_n$, $r_n = \int_{\mathbb{R}} s^n \, \mathrm{d}r(s)$, yields an expression in terms of the time $t$, the activity (mean frequency of jumps in unit time) $z$ and the 2nd, 3rd and 4th moments of the jump distribution $r$:

$$\log\left[\left(\frac{\beta}{\pi}\right)^{-d/2}\Phi_t(\phi)\right] = \frac{-\beta(tzr_2+\phi^2)}{1+\frac{\beta}{2}\frac{tzr_4+2t^2z^2r_2^2+4tzr_2\phi^2+4tzr_3\phi}{tzr_2+\phi^2}} + \mathcal{O}(\beta^3). \tag{44}$$

We have tested formula (44) by comparison with large samples of pseudo-random numbers that simulate a compound Poisson process with additional diffusion, $Z$ that is,

$$Z = X(a,\sigma^2) + s_1Y_1(z_1) - s_2Y_2(z_2), \tag{45}$$

[10] By a scaling of $\mu$ and $z$, one can always obtain $t = 1$.

[11] In fact, resummation of expansions in physics leads to results with impressive numerical precision (cf., e.g., [15], page 506). Often, validity of the resummed approximation is only undermined by phase transitions – which are, of course, absent on the discrete, finite space $\{1,\ldots,d\}$, as considered here.

where $X(a, \sigma^2)$ has Gaussian distribution with mean $a$ and standard deviation $\sigma = 1/\sqrt{2\beta}$. $Y_j(z_j)$ is Poisson distributed with parameter $z_j$, $j = 1, 2$. $s_1, s_2$ give the length of positive/negative jumps, respectively. $X, Y_1, Y_2$ are all independent. This choice is motivated by the fact that pseudo-random samples for $Z$ can be obtained using standard software routines. The parameters of the model are $z = z_1 + z_2$ and $r = (z_1\delta_{s_1} + z_2\delta_{-s_2})/(z_1 + z_2)$, and for the shift $\phi \to \phi - a$, one obtains $a = z_1 s_1 - z_2 s_2$. The normalization constant of the exponentiated right-hand side of (44) has been calculated numerically. The range of parameters has been $0.1 \lesssim \beta \lesssim 10$ and $0.5 \lesssim z \lesssim 4$, $1 \lesssim s_j \lesssim 10$, $j = 1, 2$.

If $\Psi^{\text{jump}}(-\mathrm{i}\xi)$ is symmetric under $\xi \to -\xi$ (here, $s_1 = s_2$ and $z_1 = z_2$), the results of the second order calculation are very good, as long as jumps are not too large ($\lesssim \sqrt{2\mu}$). However, the same is true for a best-fit Gaussian with mean $a$ and variance $zr_2 + 2\mu$.

The situation changes if the distribution of $Z$ becomes skew (e.g., by setting $z_2 = s_2 = 0$). In this situation, more heavy tails are being formed. A typical result is displayed in Figure 3.[12] The overall precision of the second-order calculation is again comparable with a best-fit Gaussian model, but the second-order calculation typically predicts quantiles $q_\alpha$ with $\alpha$ in the order of a few percent or above 95% with higher precision than the Gaussian model (cf. Figure 3(b)). As such quantiles are mostly used in statistics, the second order calculation already has some practical relevance.

# Appendix

## A.1. Proofs for Section 2.1

**Proof of Lemma 2.1.** $\mathrm{e}^{\mathrm{i}t\Psi(-\mathrm{i}\xi)} \to 1$ as $t \searrow 0$ in the $\mathcal{S}'(\mathbb{R}^d)$ topology, hence $\nu_t = \mathcal{F}^{-1} \times (\mathrm{e}^{t\Psi(-\mathrm{i}\xi)}) \to \mathcal{F}^{-1}(1) = \delta_0$ in $\mathcal{S}'(\mathbb{R}^d)$ for $t \searrow 0$. Here, $\delta_0$ is the Dirac measure with mass one at zero. Consequently, $\Phi_t(\phi) = \nu_t * f(\phi) \to \delta_0 * f(\phi) = f(\phi)$ as $t \searrow 0$.

It remains to show that $\Phi_t$ is differentiable in $t$ and fulfills equation (1). First, the differential quotient $\frac{\mathrm{e}^{t\Psi(-\mathrm{i}\xi)} - \mathrm{e}^{t'\Psi(-\mathrm{i}\xi)}}{t-t'}$ converges to $\Psi(-\mathrm{i}\xi)\mathrm{e}^{t\Psi(-\mathrm{i}\xi)}$ pointwise for all $\xi \in \mathbb{R}^d$. By condition (P2), we get

$$\left| \frac{\mathrm{e}^{t\Psi(-\mathrm{i}\xi)} - \mathrm{e}^{t'\Psi(-\mathrm{i}\xi)}}{t - t'} \right| \leq |\Psi(-\mathrm{i}\xi)\mathrm{e}^{t\Psi(-\mathrm{i}\xi)}|(1 + \|\xi\|)^{|t-t'|c} \tag{46}$$

and by (P2) and (P3), this implies that the differential quotient is uniformly polynomially bounded (for, say, $|t - t'| < 1$). It follows by dominated convergence that the differential quotient converges in $\mathcal{S}'(\mathbb{R}^d)$ as $t' \to t$, which implies that $\mathcal{F}(\nu_t)$ is differentiable in $t$ in

[12] Actually, in this example, $\sigma_{\text{Gauss}} = \sqrt{2\mu} = 1.58 < 12 = 2 \times$ jump length, hence the jump part for $z_1 = 2$ is not small compared with the diffusive part, which provides a nice illustration of the power of resummation.

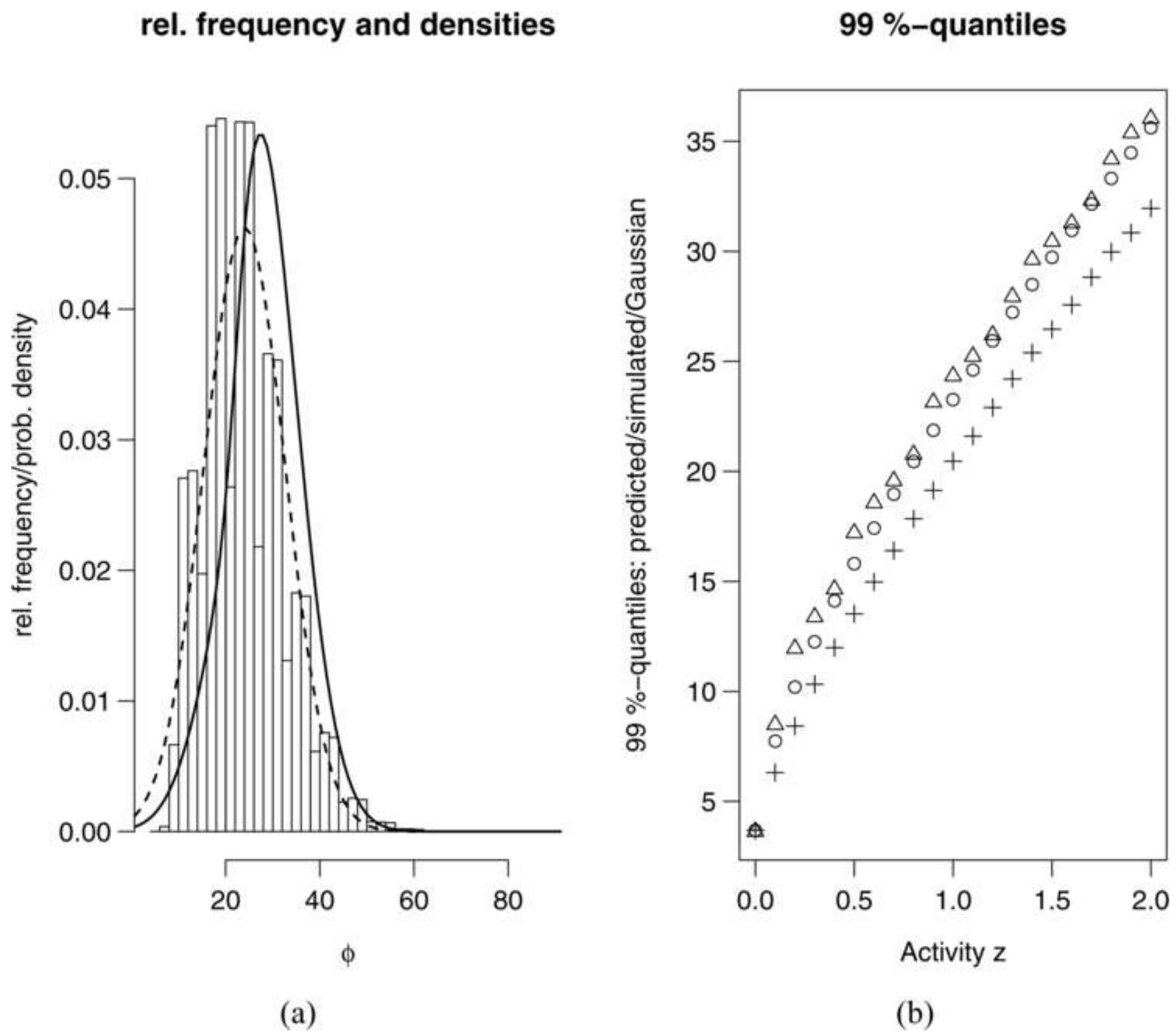

**Figure 3.** (a) Relative frequency of pseudo-random numbers (sample size $n = 10^6$), density of 2nd-order Padé approximant (44) (solid line) and a best Gaussian fit (dashed line) for $z_1 = 2, z_2 = 0, s_1 = 6$ and $\beta = 0.2$; (b) plot of 99%-sample quantiles of pseudo-random numbers (triangles), predicted 99%-quantiles from (44) (circles) and a best Gaussian fit (+). Here, the data are as in (a), but $z_1$ varies from 0 to 2.

the $\mathcal{S}'(\mathbb{R}^d)$ topology. By the continuity of $\mathcal{F}$, the same applies to $\nu_t$ and hence $\Phi_t(\phi)$ is differentiable.

It is immediate that $\mathcal{F}(\frac{\mathrm{d}}{\mathrm{d}t}\nu_t)(\xi) = \Psi(-\mathrm{i}\xi)\mathrm{e}^{t\Psi(-\mathrm{i}\xi)}$, from which we get that $\mathcal{F}(\frac{\mathrm{d}}{\mathrm{d}t}\Phi_t)(\xi) = \Psi(-\mathrm{i}\xi)\mathcal{F}(\Phi_t)(\xi)$. From the fact that the time derivative can be taken out of the Fourier transform, uniqueness of the Fourier-transformed solution follows. The argument is completed by application of the inverse Fourier transform. □

**Proof of Proposition 2.2.** In the following, let $\|h\|_{\alpha,N}$ be a Schwartz norm where $N\in\mathbb{N}$, $\alpha\in\mathbb{N}_0^d$, that is,

$$\|h\|_{\alpha,N}=\sup_{\phi\in\mathbb{R}^d}(1+\|\phi\|)^N|D_\alpha h(\phi)|. \tag{47}$$

**Lemma A.1.** *For* $\beta>0$, $\Phi_t^\beta(\phi)$ *is infinitely differentiable in* $\beta$ *and* $\frac{\mathrm{d}^m}{\mathrm{d}\beta^m}\Phi_t^\beta(\phi)=(-1)^m\nu_t*(V^mf_\beta)(\phi)$.

**Proof.** Since $\Phi_t^\beta(\phi)=\nu_t*f_\beta(\phi)$ and $\nu_t\in\mathcal{S}'(\mathbb{R}^d)$, it suffices to show that the mapping $(0,\infty)\ni\beta\to f_\beta=\mathrm{e}^{-\beta V}\in\mathcal{S}(\mathbb{R}^d)$ is $m$ times differentiable on $(0,\infty)$, $m\in\mathbb{N}$ arbitrary, and $\frac{\mathrm{d}^m}{\mathrm{d}\beta^m}f_\beta=(-1)^mV^mf_\beta$. In fact, this statement is trivial for $m=0$. Suppose that we have proven the statement up to $m-1$. We have to show that the differential quotient $(-1)^{m-1}V^{m-1}(\frac{f_\beta-f_{\beta'}}{\beta-\beta'})$ converges to $(-1)^mV^mf_\beta$ as $\beta'\to\beta$. As the common prefactor $(-1)^{m-1}V^{m-1}$ is a polynomial, it suffices to treat the case $m=1$. For $\phi\in\mathbb{R}^d$, we have

$$-V(\phi)f_\beta(\phi)-\frac{f_\beta(\phi)-f_{\beta'}(\phi)}{\beta-\beta'}=-V(\phi)\frac{\int_\beta^{\beta'}(f_\beta-f_\tau)\,\mathrm{d}\tau}{\beta-\beta'}. \tag{48}$$

Hence,

$$\begin{aligned}&\left\|-Vf_\beta(\phi)-\frac{f_\beta-f_{\beta'}}{\beta-\beta'}\right\|_{\alpha,N}\\&\quad\le K\sum_{\alpha'\in\mathbb{N}_0^d:|\alpha'|<|\alpha|}\sup_{\phi\in\mathbb{R}^d}(1+\|\phi\|)^{N_{\alpha'}}\sup_{\tau\in(\beta,\beta')}|D_{\alpha'}f_\beta(\phi)-D_{\alpha'}f_\tau(\phi)|\end{aligned} \tag{49}$$

for $K$ and $N_{\alpha'}$ sufficiently large. The derivative $D_{\alpha'}f_\beta(\phi)$ is of the form $P_{\alpha'}(\phi,\beta)f_\beta(\phi)$, where $P_{\alpha'}(\phi,\beta)$ is a polynomial in $\beta$ and $\phi$. Hence, for $|\beta-\beta'|<1$, $\beta,\beta'>\epsilon>0$,

$$\begin{aligned}&\sup_{\tau\in(\beta,\beta')}|D_{\alpha'}f_\beta(\phi)-D_{\alpha'}f_\tau(\phi)|\\&\quad=\sup_{\tau\in(\beta,\beta')}|P_{\alpha'}(\phi,\beta)f_\beta(\phi)-P_{\alpha'}(\phi,\tau)f_\tau(\phi)|\\&\quad\le\sup_{\tau\in(\beta,\beta')}|P_{\alpha'}(\phi,\beta)-P_{\alpha'}(\phi,\tau)|f_\beta(\phi)+\sup_{\tau\in(\beta,\beta')}|P_{\alpha'}(\phi,\tau)(f_\beta(\phi)-f_\tau(\phi))'|\\&\quad\le K_{\alpha'}|\beta-\beta'|(1+\|\phi\|)^{Q_{\alpha'}}(f_\beta(\phi)+|V(\phi)|\mathrm{e}^{(\beta+1)L}\mathrm{e}^{-\epsilon\langle\lambda^{(\bar p)},\phi^{\otimes\bar p}\rangle_{\bar p}}),\end{aligned} \tag{50}$$

with $L=\min_{\phi\in\mathbb{R}^d}V(\phi)-\epsilon\langle\lambda^{(\bar p)},\phi^{\otimes\bar p}\rangle_{\bar p}$, $K_{\alpha'},Q_{\alpha'}$ sufficiently large. Multiplying by $(1+\|\phi\|)^{N_{\alpha'}}$ and taking the supremum over $\phi\in\mathbb{R}^d$ implies that the right-hand side of (49) vanishes as $\beta'\to\beta$. This completes the proof. □

It remains to show that the derivatives of $\Phi_t^\beta(\phi)$ extend continuously to $\beta=0$. We start with two technical lemmas.

Let $h_r \in \mathcal{S}(\mathbb{R}^d)$, $r \in (0,\infty)$, be a family of functions such that $h_r \to h$ for $r \searrow 0$ in $\mathcal{S}(\mathbb{R}^d)$, where $\int_{\mathbb{R}^d} h\,\mathrm{d}\xi = 1$. Further, let $h^\epsilon(\xi) = \epsilon^{-d} h(\xi/\epsilon)$ and $h_r^\epsilon(\xi)$ be defined likewise. Let $r(\epsilon)$ be a function such that $r(\epsilon) \to 0$ as $\epsilon \searrow 0$.

**Lemma A.2.** *Let $F : \mathbb{R}^d \to \mathbb{C}$ be a function that is polynomially bounded and has polynomially bounded partial derivatives. Let $h_r^\epsilon$ and $r(\epsilon)$ be defined as above. Then $\lim_{\epsilon\searrow 0} \int_{\mathbb{R}^d} F(\xi) \times h_{r(\epsilon)}^\epsilon(\xi)\,\mathrm{d}\xi = F(0)$.*

**Proof.** We consider the following estimate:

$$\begin{aligned} &\left| \int_{\mathbb{R}^d} F(\xi) h_{r(\epsilon)}^\epsilon(\xi)\,\mathrm{d}\xi - F(0) \right| \\ &\quad = \left| \int_{\mathbb{R}^d} F(\epsilon\xi) h_{r(\epsilon)}(\xi)\,\mathrm{d}\xi - F(0) \right| \\ &\quad \le \left| \int_{\mathbb{R}^d} F(\epsilon\xi)(h_{r(\epsilon)}(\xi) - h(\xi))\,\mathrm{d}\xi \right| + \left| \int_{\mathbb{R}^d} (F(\epsilon\xi) - F(0)) h(\xi)\,\mathrm{d}\xi \right|. \end{aligned} \tag{51}$$

The first term on the right-hand side can be estimated as follows:

$$\left| \int_{\mathbb{R}^d} F(\epsilon\xi)(h_{r(\epsilon)}(\xi) - h(\xi))\,\mathrm{d}\xi \right| \le \int_{\mathbb{R}^d} \frac{|F(\epsilon\xi)|}{(1+\|\xi\|)^N}\,\mathrm{d}\xi \|h_{r(\epsilon)} - h\|_{0,N}. \tag{52}$$

Recall that $F$ is polynomially bounded. Thus, for $N$ sufficiently large, the integral on the right-hand side converges to $|F(0)| \int_{\mathbb{R}^d} \frac{1}{(1+\|\xi\|)^N}\,\mathrm{d}\xi$ by dominated convergence. Hence, the right-hand side of (52) tends to zero as $\epsilon \searrow 0$.

The second term on the right-hand side of (51), for $0 < \epsilon < 1$, can be estimated by

$$\begin{aligned} \left| \int_{\mathbb{R}^d} (F(\epsilon\xi) - F(0)) h(\xi)\,\mathrm{d}\xi \right| &= \left| \int_0^\epsilon \int_{\mathbb{R}^d} \langle \nabla F(s\xi), \xi \rangle h(\xi)\,\mathrm{d}\xi\,\mathrm{d}s \right| \\ &\le \epsilon K \int_{\mathbb{R}^d} (1+\|\xi\|)^N |h(\xi)|\,\mathrm{d}\xi, \end{aligned} \tag{53}$$

for $K$ and $N$ sufficiently large, by the polynomial boundedness of the partial derivatives of $F$. Hence, the second term on the right-hand side of (51) also tends to zero as $\epsilon \searrow 0$. □

**Lemma A.3.** *Let $V_\epsilon(\varphi) = \langle \lambda^{(\bar p)}, \varphi^{\otimes \bar p} \rangle_{\bar p} + \sum_{p=0}^{\bar p - 1} \langle \lambda_\epsilon^{(p)}, \varphi^{\otimes p} \rangle_p$ such that $\lambda_\epsilon^{(p)} \to 0$ as $\epsilon \searrow 0$, $p = 0, \ldots, \bar p - 1$ and $f^\epsilon = \mathrm{e}^{-V_\epsilon}$. Let $f^0(\varphi) = \mathrm{e}^{-\langle \lambda^{(\bar p)}, \varphi^{\otimes \bar p} \rangle_{\bar p}}$. Then $f^\epsilon \to f^0$ in $\mathcal{S}(\mathbb{R}^d)$ as $\epsilon \searrow 0$.*

**Proof.** We first assume $\alpha = 0$, then $|f^0 - f^\epsilon| \le \int_0^1 |Q\mathrm{e}^{-sQ}|\,\mathrm{d}s\, f^0 \le |Q|\mathrm{e}^{|Q|} f^0$, with $Q$ a polynomial that has vanishing coefficients as $\epsilon \searrow 0$. Multiplying this by another polynomial and taking the supremum over $\mathbb{R}^d$ thus gives an expression that vanishes in that limit. Hence, $\|f^\epsilon - f^0\|_{0,N} \to 0$ for $\epsilon \searrow 0$.

Next, let $\alpha\in\mathbb{N}_0^d$, $\alpha\neq 0$. We note that each contribution to $D_\alpha(f_0-f^\epsilon)$ is of the form $f_0\mathrm{e}^Q P$ or of the form $S(f^0-f^\epsilon)$. $Q$ is of degree at most $\bar p-1$ (potentially, $Q=0$). $P,Q$ and $S$ are polynomials. The coefficients of $P$ and $Q$ are vanishing as $\epsilon\searrow 0$. Hence, we can conclude, as above, that $\|f^\epsilon-f^0\|_{\alpha,N}\to 0$ as $\epsilon\searrow 0$. □

We now conclude the proof of Proposition 2.2. By Plancherel's formula, (P2)–(P3) and Lemma A.1, we have

$$\frac{\mathrm{d}^m}{\mathrm{d}\beta^m}\Phi_t^\beta(\phi)=(-1)^m\nu_t*(V^m f_\beta)(\phi)=(-1)^m\int_{\mathbb{R}^d}\mathcal{F}^{-1}(f_{\beta,\phi})(\xi)V_\phi^m(\mathrm{i}\nabla_\xi)\mathrm{e}^{t\Psi(-\mathrm{i}\xi)}\,\mathrm{d}\xi. \quad (54)$$

For a polynomial $W(\varphi)=\sum_{p=0}^{\bar p}\langle w^{(p)},\varphi^{\otimes p}\rangle_p$, we set $W_\beta(\varphi)=\sum_{p=0}^{\bar p}\beta^{1-p/\bar p}\langle w(p),\varphi^{\otimes p}\rangle_p$. Let $f^0(\varphi)=\mathrm{e}^{-\langle\lambda^{(\bar p)},\varphi^{\otimes\bar p}\rangle_{\bar p}}$ and let $f^\beta(\phi)=\mathrm{e}^{-(V_\phi)_\beta}$. By Lemma A.3, $f^\beta\to f^0$ in $\mathcal{S}(\mathbb{R}^d)$ as $\beta\searrow 0$. Hence, $\mathcal{F}^{-1}(f_\phi^\beta)\to\mathcal{F}^{-1}(f^0)$ in $\mathcal{S}(\mathbb{R}^d)$ and $\int_{\mathbb{R}^d}\mathcal{F}^{-1}(f^0)\,\mathrm{d}\xi=f^0(0)=1$.

Furthermore, $f_{\beta,\phi}(\varphi)=\mathrm{e}^{-\beta V_\phi(\varphi)}=\mathrm{e}^{-(V_\phi)_\beta(\beta^{1/\bar p}\varphi)}=f_\phi^\beta(\beta^{1/\bar p}\varphi)$. Inserting this into the right-hand side of (54), we obtain

$$\frac{\mathrm{d}^m}{\mathrm{d}\beta^m}\Phi_t^\beta(\phi)=(-1)^m\int_{\mathbb{R}^d}\beta^{-d/\bar p}\mathcal{F}^{-1}(f_\phi^\beta)(\xi/\beta^{1/\bar p})V_\phi^m(\mathrm{i}\nabla_\xi)\mathrm{e}^{t\Psi(-\mathrm{i}\xi)}\,\mathrm{d}\xi. \quad (55)$$

Note that by (P2)–(P3), $V_\phi(\mathrm{i}\nabla_\xi)\mathrm{e}^{t\Psi(-\mathrm{i}\xi)}$ as a function in $\xi$ is polynomially bounded with polynomially bounded derivatives. We can thus apply Lemma A.2 and we get

$$\lim_{\beta\searrow 0}\frac{\mathrm{d}^m}{\mathrm{d}\beta^m}\Phi_t^\beta(\phi)=(-1)^m V_\phi^m(\mathrm{i}\nabla_\xi)\mathrm{e}^{t\Psi(-\mathrm{i}\xi)}|_{\xi=0}=(-1)^m\langle V_\phi^m(\varphi)\rangle_{\nu_t}. \quad (56)$$

This completes the proof. □

## A.2. Notions from graph theory

For the reader's convenience, we develop some graph-theoretic notions, some of them being non-standard, essentially following [6].

Given a set $M$, let $M_s^2$ be the set of subsets of $M$ containing two elements. A *graph* $G$ over $M$ is then a subset of $M_s^2$. The points of $M$ are called *vertices* of $G$ and the elements (i.e., non-ordered pairs of $M$) $\{m_1,m_2\}\in G$ are called the *edges* (or bonds). $\mathcal{G}(M)$ stands for the collection of all graphs over $M$. We could also replace $M_s^2$ with $M^2=M\times M$, the set of ordered pairs. We then obtain a graph with *directed edges*.

Let $N\subseteq M$. The following construction explains the notion of a *graph with distinguishable legs* (endpoints of edges) at the vertices in $N$. Let $\hat M$ be a set, $\hat\pi:\hat M\to M$ a surjective map and $\hat\pi^2:\mathcal{G}(\hat M)\to\mathcal{G}(M)$ the induced mapping on graph level defined by application of $\hat\pi$ on the endpoints of edges. Let $G\in\mathcal{G}(M)$ be fixed and for $m\in M$, let $p(m,G)$ be the number of edges connected with $m$. Let $(\hat M_G,\hat\pi_G)$ be as above such that for $m\in M\setminus N$, $\sharp\hat\pi_G^{-1}(\{m\})=1$ and $\sharp\hat\pi_G^{-1}(\{m\})=p(m,G)$ for $m\in N$. $\sharp A$ denotes the cardinality of the set $A$. A *leg* $\hat{m}$ at a vertex $m$ is a point in $\hat{\pi}_G^{-1}(\{m\})$. The class of graphs with distinguishable legs at vertices in $N$ associated with an ordinary graph $G \in \mathcal{G}(M)$ is then given by $[G|N] = \{\hat{G} \in (\hat{\pi}_G^2)^{-1}(G) : \forall m \in N, \hat{m} \in \hat{\pi}^{-1}(m), \sharp\{e \in \hat{G} : \hat{m} \in e\} = 1\}$. The class of graphs with distinguishable legs at the vertices in $N$ is then $[\mathcal{G}|N](M) = \bigcup_{G \in \mathcal{G}(M)} [G|N]$.

Let $L \subseteq M$ and $s$ be in the permutation group of $L$, $\mathrm{Perm}(L)$. Obviously, $s$ extends to a permutation on $M$ by putting $s(m) = m$ for $m \in M \setminus L$. Let $s^2 : \mathcal{G}(M) \to \mathcal{G}(M)$ be the induced mapping on graph level. We can then consider the equivalence classes under exchange of vertices in $L$, $[\mathcal{G} : L](M) = \mathcal{G}(M)/(\mathrm{Perm}(L))^2$. Such an equivalence class is called a *graph with indistinguishable vertices* in $L$. It is easy to check that the above two constructions are compatible, that is, a vertex can belong to a set $L$ of indistinguishable vertices with distinguishable or indistinguishable legs and a distinguishable vertex can have distinguishable or indistinguishable legs.

A *graph with several types of vertices* has sets of vertices $M_1, \dots, M_k$ and is a graph over $M = \dot{\bigcup}_{l=1}^k M_k$.

Finally, a *topological graph* is a graph with indistinguishable vertices in each class of vertices and with indistinguishable legs. The *multiplicity* of a topological graph, given some other class of graphs (generalized Feynman graphs in the present case), is the number of such graphs that coincide with the topological graph up to the labeling of indistinguishable vertices and legs.

Let us also note that graph-theoretic notions like the *connectedness* of a graph with distinguishable legs $\hat{G}$ are meant in the sense that the associated "projected" graph $G$ or even topological graph $G_{\mathrm{top}}$ are connected.[13]

## A.3. Padé resummation

The idea of Padé resummation is to replace the Taylor polynomial with a rational function that has the given Taylor approximation at the origin but often has a better large $\beta$ behavior. We essentially follow [10], Chapter 16.2. The $[M/N]$-Padé approximant is defined as

$$[M/N](\beta) = \frac{a_0 + a_1\beta + \cdots + a_M\beta^M}{1 + b_1\beta + \cdots + b_N\beta^N}. \tag{57}$$

The $M + N + 1$ coefficients $a_0, \dots, a_M; b_1, \dots, b_N$ are determined by the equation

$$[M/N](\beta) = \sum_{m=0}^{M+N} h_m\beta^m + \mathcal{O}(\beta^{N+M+1}). \tag{58}$$

[13] $\hat{G}$ as a graph in $\mathcal{G}(\hat{M})$ might be disconnected even if $G$ is connected.

Multiplication with the denominator of $[M/N]$ and comparison of coefficients leads to the system of equations

$$\left.\begin{array}{l} b_N h_{M-N+1} + b_{N-1} h_{M-N+2} + \cdots + b_0 h_{M+1} = 0 \\ \quad\vdots \qquad\qquad \vdots \qquad\qquad \vdots \quad \vdots \\ b_N h_M + b_{N-1} h_{M+1} + \cdots + b_0 h_{M+N} = 0 \end{array}\right\} b_0 = 1 \tag{59}$$

and

$$\left.\begin{array}{r} b_0 h_0 = a_0 \\ b_1 h_0 + b_0 h_1 = a_1 \\ \vdots \\ b_M h_0 + b_{M-1} h_1 + \cdots + b_0 h_M = a_M \end{array}\right\} b_0 = 1. \tag{60}$$

Solving for $N = M = 1$, as required in Section 3.2, one obtains

$$\begin{aligned} b_1 &= -\frac{h_2}{h_1}, \\ a_0 &= h_0, \\ a_1 &= h_1 - \frac{h_0 h_2}{h_1}. \end{aligned} \tag{61}$$

With $h_0, h_1, h_2$ as in (42), this leads to

$b_1 = \frac{1}{2}$ 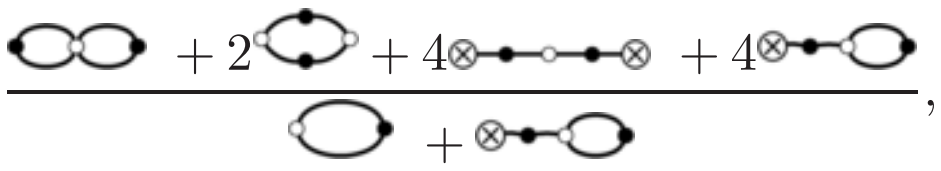,

$a_0 = 0,$ (62)

$a_1 = -[$ 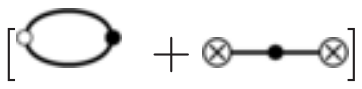 $],$

from which (43) follows.

# Acknowledgements

During his stay in Bonn, B. Smii received financial support from the DAAD. When this article was prepared, H. G. was supported by the D.F.G. project "Stochastic methods in quantum field theory." H. T. received financial support through the SFB 611. The authors gratefully acknowledge this support.

The authors also thank Sergio Albeverio for his interest and encouragement while this work was underway. Habib Ouerdiane contributed interesting discussions and comments on the topic during his stay in Bonn in the summer of 2005.